\documentclass[10pt]{amsart}
\newcommand{\proclaim}[2]{\medbreak {\bf #1}{\sl #2} \medbreak}

\newcommand{\ntop}[2]{\genfrac{}{}{0pt}{1}{#1}{#2}}

\let\newpf\proof \let\proof\relax 
\newenvironment{pf}{\newpf[\proofname]}{\qed\endtrivlist}

\def\d{{\underline d}}

\def\n{{\bf n}}

\def\g{{\gamma}}

\def\NR{{\mathcal {NR}}}

\def\UF{{\mathbb {UF}}}

\def\O{{\mathbb O}}

\def\op{{\overline \partial}}

\def\0{{\mathbf 0}}

\def\hh{{\mathbf h}}

\def\cal{\mathcal}

\newtheorem{thm}{Theorem}[section]

\newtheorem{lemma}[thm]{Lemma}

\newtheorem{prop}[thm]{Proposition}

\theoremstyle{remark}
\newtheorem{rem}{Remark}[section]

\numberwithin{equation}{section}

\def\conj{\operatorname{conj}}

\def \bn {\hfill \\ \smallskip\noindent}

\theoremstyle{definition}

\newtheorem{definition}{Definition}[section]
\def\proof{\bn {\bf Proof.} }

\def\note#1
{\marginpar
{\tiny $\leftarrow$
\par
\hfuzz=20pt \hbadness=9000 \hyphenpenalty=-100 \exhyphenpenalty=-100
\pretolerance=-1 \tolerance=9999 \doublehyphendemerits=-100000
\finalhyphendemerits=-100000 \baselineskip=6pt
#1}\hfuzz=1pt}

\newcommand{\ra}{\rightarrow}

\newcommand{\inter}{\operatorname{int}}
\renewcommand{\mod}{\operatorname{mod}}

\newcommand{\orb}{\operatorname{orb}}

\newcommand{\id}{\operatorname{id}}

\newcommand{\Dil}{\operatorname{Dil}}

\newcommand{\Om}{{\Omega}}

\newcommand{\AAA}{{\cal A}}

\newcommand{\EE}{{\cal E}}

\newcommand{\FF}{{\cal F}}

\newcommand{\LL}{{\cal L}}
\newcommand{\MM}{{\cal M}}

\newcommand{\RR}{{\cal R}}
\newcommand{\TT}{{\cal T}}

\newcommand{\UU}{{\cal U}}
\newcommand{\VV}{{\cal V}}
\newcommand{\WW}{{\cal W}}
\newcommand{\XX}{{\cal X}}
\newcommand{\YY}{{\cal Y}}

\newcommand{\C}{{\mathbb C}}
\newcommand{\D}{{\mathbb D}}
\newcommand{\E}{{\mathbb E}}
\newcommand{\F}{{\mathbb F}}

\newcommand{\R}{{\mathbb R}}

\newcommand{\U}{{\mathbb U}}

\newcommand{\Z}{{\mathbb Z}}

\def\B0{{\bold{0}}}

\catcode`\@=12

\def\Empty{}
\newcommand\oplabel[1]{
  \def\OpArg{#1} \ifx \OpArg\Empty {} \else
  	\label{#1}
  \fi}
		
\newcommand{\comm}[1]{}
\newcommand{\comment}[1]{}

\begin{document}

\title[Phase-Parameter relation]
{Phase-parameter relation and sharp statistical properties for general
families of unimodal maps}

\author{Artur Avila and Carlos Gustavo Moreira}

\address{
Coll\`ege de France -- 3 Rue d'Ulm \\
75005 Paris -- France.
}
\email{avila@impa.br}

\address{
IMPA -- Estr. D. Castorina 110 \\
22460-320 Rio de Janeiro -- Brazil.
}
\email{gugu@impa.br}

\thanks{Partially supported by Faperj and CNPq, Brazil.}   

\begin{abstract}

We obtain estimates relating the phase space and the parameter space
of analytic families of unimodal maps.
Using those estimates, we show that typical analytic unimodal maps admit
a quasiquadratic renormalization.  This reduces the study of the statistical
properties of typical unimodal maps to the quasiquadratic case which had
been studied in \cite {AM2}.  The estimates proved here correspond exactly
to the Phase-Parameter relation proved in \cite {AM} in the quadratic case,
and allows one to obtain sharp estimates on the dynamics of typical
unimodal maps which were available only
in the quadratic case: as an example we
conclude that the exponent of the polynomial recurrence of the critical
orbit is exactly one.  We also show that those ideas lead
to a new proof of a Theorem of Shishikura: the set of
non-renormalizable parameters in the boundary of the Mandelbrot set has
Lebesgue measure zero.  Further applications of those results can be found
in the companion paper \cite {AM3}.

\end{abstract}

\setcounter{tocdepth}{1}

\date\today

\maketitle

\tableofcontents

\section{Introduction}

A unimodal map is a smooth (at least $C^2$) map $f:I \to I$, where $I
\subset \R$ is an interval, which has one unique critical point
$c \in \inter I$ which is a maximum. 
Let us say that $f$ is regular if it has a quadratic critical point, is
hyperbolic and its critical point is not periodic or preperiodic.
By a result of Kozlovski \cite {K2},
the set of regular maps coincide with the set of
structurally stable unimodal maps, and it follows that
the set of regular maps is open and dense in all smooth (and even analytic)
topologies.

The most studied family of unimodal maps is the quadratic family
$p_\lambda=\lambda-x^2$, $-1/4 \leq \lambda \leq 2$.
In \cite {AM} it was shown that for a typical non-regular quadratic map
$p_{\lambda_0}$, the phase space of $p_{\lambda_0}$ near the critical point
$0$ and the parameter space near $\lambda_0$ are related by
some metric rules called the {\it Phase-Parameter relation}
(notice that it is crucial that the phase and parameter of
the quadratic family have the same dimension).
The proof of \cite {AM} was tied to the combinatorial theory of the
Mandelbrot set, so it can only work for quadratic maps
(or, at most, full unfolded families of quadratic-like maps, see \cite
{parapuzzle}).

Let us say that an analytic family of unimodal maps is non-trivial if
regular parameters are dense (in particular non-trivial analytic
families are dense in any topology).
The first main result of this paper is the following (see \S \ref
{phase-parameter} for the precise definition of the Phase-Parameter
relation):

\proclaim{Theorem A.}
{
Let $f_\lambda$ be a one-parameter non-trivial analytic
family of unimodal maps.  Then $f_\lambda$
satisfies the Phase-Parameter relation at
almost every non-regular parameter.
}

The Phase-Parameter relation has many remarkable consequences for the study
of the dynamical behavior of typical parameters.
Our second main result is an application of the Phase-Parameter relation:

\proclaim{Theorem B.}
{
Let $f_\lambda$ be a non-trivial analytic
family of unimodal maps (any number of parameters).  Then almost every
parameter is either regular or has a renormalization which is topologically
conjugate to a quadratic polynomial.
}

This result allows one to reduce the study of typical unimodal maps to the
special case of unimodal maps which are quasiquadratic
(persistently topologically conjugate to a quadratic polynomial).

\subsection{Application: statistical properties of typical unimodal maps}

Typical quasiquadratic maps had been previously studied in \cite {AM2}.
Their main result is that the dynamics of
typical quasiquadratic maps have an excellent statistical
description (in terms of physical
measures, decay of correlations and stochastic stability), thus
answering the Palis Conjecture (see \cite {AM2} for details)
in the unimodal quasiquadratic case.

For regular
maps, the good statistical description comes for free.
For a non-regular map $f$,
it is related to essentially two properties regarding its
critical point $c$: the Collet-Eckmann
condition\footnote{A unimodal map $f$ is Collet-Eckmann if $|Df^n(f(c))|>C
\lambda^n$ for some constants $C>0$ and $\lambda>1$.} and subexponential
recurrence\footnote{That is, for every $\alpha>0$, $|f^n(c)-c|>e^{-\alpha
n}$ for $n$ sufficiently big.}.

Thus, \cite {AM2} achieves the good statistical description via a dichotomy:
typical quasiquadratic maps are either regular or Collet-Eckmann and
subexponentially recurrent.  This is done in both the analytic case and
the smooth case ($C^k$, $k=3,...,\infty$).  For typical non-regular
analytic quasiquadratic maps, it is proved even more, that the critical point
is polynomially recurrent\footnote{That is, there exists
$\gamma>0$ such that $|f^n(c)-c|>n^{-\gamma}$ for every
$n$ sufficiently big.}.

Our Theorem B allows us to immediately obtain the analytic case
in our more general setting (see Theorem \ref {C precise} for a more precise
statement):

\proclaim{Corollary C.}
{
Let $f_\lambda$ be a non-trivial analytic family of
unimodal maps (in any number of parameters).
Then almost every non-regular parameter is Collet-Eckmann and its critical
point is polynomially recurrent.
}

This allows us not only to generalize the smooth case of \cite {AM2} besides
quasiquadratic maps, but to reduce the
differentiability requirements, including the $C^2$ case in the description
(see Theorem \ref {D precise} for a more precise statement):

\proclaim{Corollary D.}
{
In generic smooth ($C^k$, $k=2,...,\infty$) families of unimodal maps (any
number of parameters),
almost every parameter is regular, or has a renormalization which is
conjugate to a quadratic map, is Collet-Eckmann and its critical point is
subexponentially recurrent.
}


\begin{rem}

The dichotomies in Corollaries C and D imply that the
dynamics of typical non-regular unimodal maps have the same
excellent statistical description of the quasiquadratic case studied by
\cite {AM2}, see also Remark \ref {ergodic} for a list of references.  In
particular, our Corollaries C and D give an answer to the Palis Conjecture
in the general unimodal case.

\end{rem}

\subsection{Sharpness}

The Phase-Parameter relation allows one to obtain very
precise estimates on the
dynamics of typical parameters.  For instance, the statistical
analysis of \cite {AM} could compute the exact exponent of the
polynomial recurrence\footnote{The exponent of the polynomial recurrence of
the critical point $c$ of a
unimodal map $f$ is the infimum of all $\g>0$ such that, for
$n$ sufficiently big, $|f^n(c)-c|>n^{-\g}$.} in the case of
the quadratic family.  The method used in \cite {AM2} to
extend results from the quadratic family to other non-trivial
families of quasiquadratic maps
(based on comparison of the respective parameter spaces)
introduces unavoidable distortion and can not be used to estimate the
exponent of the recurrence even in the quasiquadratic case.
Our Theorem A implies that the same sharp estimates obtained for the
quadratic family remain valid in general.

\proclaim{Corollary E.}
{
Let $f_\lambda$ be a non-trivial analytic
family of unimodal maps (any number of parameters).  Then almost every
parameter is either regular or has a polynomially
recurrent critical point with exponent 1.
}

We call the attention of the reader to the companion paper \cite {AM3}
where much more refined statistical applications of Theorem A are obtained.
Those results are inaccessible by the methods of \cite {AM2}, and indeed
are used to show the limitations of estimates based on comparison of
parameter spaces of different families with respect to direct
Phase-Parameter estimates.

\subsection{Complex parameters}

A very natural question raised by the description of typical parameters in
the real quadratic family is if the results generalize to complex
parameters.  It is widely expected that the description should be actually
simpler: almost every complex parameter should be hyperbolic.  However, only
partial results are available.

In this direction, let us remark that the argument of the
proof of Theorem B can be also applied in the complex setting,
and leads to a new proof of the following result of Shishikura
(unpublished, a sketch can be found as Theorem 4 in \cite {Sh}):

\proclaim{Theorem F.}
{
The set of non-hyperbolic, non-infinitely renormalizable complex quadratic
parameters has zero Lebesgue measure.
}

We discuss this application in Appendix B.

\subsection{Outline of the proof of Theorem A}

The proof of Theorem A can be divided in four parts.  The crucial step of
this paper is step (2) below, which allows us to integrate the work of
\cite {AM} and \cite {ALM}.

\noindent (1)\, Following \cite {parapuzzle} and the Appendix of \cite {AM},
we describe a complex analogous of the
Phase-Parameter relation for certain families
of complex return type maps, which model complex extensions of
the return maps $R_n:I_n \to I_n$ to the principal nest of a unimodal map
$f$.  This study is restricted to the class of so called full families.

\noindent (2)\, We show that through any given analytic unimodal map $f$
(assumed finitely renormalizable with a recurrent critical point),
there exists an analytic family $\tilde f_\lambda$ (constructed explicitly)
which gives rise (after a generalized renormalization procedure)
to a full family of complex return type maps.
Using step (1), we conclude that the
Phase-Parameter relation is valid at $f$ for this special family
$\tilde f_\lambda$.
By construction, this family is tangent to a certain special
infinitesimal perturbation considered in \cite {ALM}, where
this perturbation had been shown
to be transverse to the topological class of $f$ (which is a
codimension-one analytic submanifold).

\noindent (3)\, We show that if the Phase-Parameter relation is valid for
one transverse family at $f$, then it is valid for all transverse families at
$f$.  This step is heavily based on the results of \cite {ALM}: in order to
compare the parameter space of both families, one uses the local holonomy of
the lamination associated to the partition of spaces of unimodal maps in
topological classes.

\noindent (4)\, Using a simple generalization of
\cite {ALM} we conclude that a non-trivial family of unimodal maps is
transverse to the topological class of almost every non-regular parameter,
and that
typical parameters are finitely renormalizable with a recurrent critical
point.  This concludes the proof of Theorem A.

\subsection{Structure of the paper}

In \S \ref {preliminaries}
we give some basic background on quasiconformal maps and holomorphic
motions.  In \S \ref {compl dyn}, we discuss the
dynamics of families of complex return-type maps
(this is based on \cite {parapuzzle}) and obtain some Phase-Parameter
estimates in this context (following the sketch of the Appendix
of \cite {AM}).
In \S \ref {puzzle}
we present the results of Lyubich in \cite {puzzle} and \cite {parapuzzle}
in the generality needed for our applications.
In \S \ref {basic} we present the
basic theory of unimodal maps, and in \S \ref {hybrid}
we introduce the results of \cite {ALM} on the lamination structure
of topological classes of unimodal maps and state some
straightforward generalizations (some details are given in Appendix A).
In \S \ref {family}
we construct a special analytic family of unimodal maps which induce
a full family of return type maps and in \S \ref {phase-parameter}
we state and prove the Phase-Parameter relation for the special family.
In \S \ref {thmA} and \S \ref {thmB} we prove
Theorems A and B, and in \S \ref {pf} we show the
relation to the corollaries.
In Appendix B we give a proof of Theorem F.

{\bf Acknowledgements:}
Most of the results of this paper were announced in \cite {Av2}, and,
together with \cite {AM}, \cite {AM2} and \cite {ALM}, formed the thesis of
the first author.  The first author would like to thank
Welington de Melo and Mikhail Lyubich who collaborated in those works,
and to Viviane Baladi and Jean-Christophe Yoccoz for useful conversations.

\section{Preliminaries} \label {preliminaries}

\subsection{General notation}

Let $\Omega$ be the set of finite sequences (possibly empty)
of non-zero integers $\d=(j_1,...,j_m)$.

A {\it Jordan curve} $T$ is a subset of $\C$ homeomorphic to a circle.
A {\it Jordan disk} is a bounded open subset
$U$ of $\C$ such that $\partial U$ is a Jordan curve.

We let $\D_r(w)=\{z \in \C||z-w|<r\}$.  Let $\D_r=\D_r(0)$, and $\D=\D_1$.
If $r>1$, let $A_r=\{z \in \C|1<|z|<r\}$.
An {\it annulus} $A$ is a subset of
$\C$ such that there exists a conformal map from
$A$ to some $A_r$.  In this case, $r$ is uniquely defined
and we denote the {\it modulus} of $A$ as $\mod (A)=\ln(r)$.

\comm{
If $U$ is a Jordan disk, we say that $U$ is the {\it filling}
of $\partial U$.  A subset $X \subset U$ is said to be {\it bounded}
by $\partial U$.
}

\subsection{Graphs and sections}

Let us fix a complex Banach space $\E$.
If $\Lambda \subset \E$, a {\it graph} of a continuous map
$\phi:\Lambda \to \C$ is the set of all $(z,\phi(z)) \in \E \oplus \C$,
$z \in \Lambda$.

Let $\0:\E \to \E \oplus \C$ be defined by $\0(z)=(z,0)$.

Let $\pi_1:\E \oplus \C \to \E$, $\pi_2:\E \oplus \C \to \C$
be the coordinate projections.  Given a set
$\XX \subset \E \oplus \C$ we denote its fibers
$X[z]=\pi_2(\XX \cap \pi_1^{-1}(z))$.

A {\it fiberwise map} $\FF:\XX \to \E \oplus \C$ is a map such that
$\pi_1 \circ \FF=\pi_1$.  We denote its fibers
$F[z]:X[z] \to \C$ such that $\FF(z,w)=(z,F[z](w))$.

Let $B_r(\E)$ be the ball of radius $r$ around $0$.

\subsection{Quasiconformal and quasisymmetric maps}
\label {quasiconformal maps}

Let $U \subset \C$ be a domain.  A map $h:U \to \C$ is
{\it $K$-quasiconformal} ($K$-qc) if it is a homeomorphism onto its image
and for any annulus $A \subset U$,
$\mod (A)/K \leq \mod (h(A)) \leq K \mod (A)$.
The minimum such $K$ is called the {\it dilatation} $\Dil(h)$ of $h$.

\comm{
Let $h:X \to \C$ be a homeomorphism and let $C,\epsilon>0$.  An extension
$H:U \to \C$ of $h$ to a Jordan disk $U$ is $(C,\epsilon)$-qc
if there exists an annulus $A \subset U$ with $\mod (A)>C$
such that $X$ is contained in the bounded component of the complement of
$A$ and $H$ is $1+\epsilon$-qc.
}

A homeomorphism $h:\R \to \R$ is said to be $\g$-quasisymmetric if it has a
real-symmetric extension $h:\C \to \C$ which is quasiconformal with
dilatation bounded by $\g$.  If $X \subset \R$, we will also say that
$h:X \to \R$ is $\g$-qs if it has a $\g$-qs extension.

A quasiconformal vector field $\alpha$ of $\C$ is a continuous
vector field with locally integrable distributional derivatives
$\op \alpha$ and $\partial \alpha$
in $L^1$ and $\op \alpha \in L^\infty$.

\subsection{Holomorphic motions} \label {holomorphic motions}

Let $\Lambda$ be a connected open set of a Banach space $\E$.
A {\it holomorphic motion} $h$ over $\Lambda$ is a family of
holomorphic maps defined on $\Lambda$ whose graphs
(called {\it leaves} of $h$) do not intersect.
The {\it support} of $h$
is the set $\XX \subset \C^2$ which is the union of the leaves of $h$.

The {\it transition} (or holonomy) maps
$h[z,w]:X[z] \to X[w]$, $z,w \in \Lambda$,
are defined by $h[z,w](x)=y$ if $(z,x)$ and $(w,y)$ belong to the same leaf.

Given a holomorphic motion $h$ over a domain $\Lambda$, a holomorphic motion
$h'$ over a domain $\Lambda' \subset \Lambda$ whose leaves are contained in
leaves of $h$ is called a {\it restriction} of $h$.  If $h$ is a
restriction of $h'$ we also say that $h'$ is an {\it extension} of $h$.

Let $K:[0,1) \to \R$ be defined by $K(r)=(1+\rho)/(1-\rho)$ where $0 \leq
\rho <1$ is such that the hyperbolic distance in $\D$ between
$0$ and $\rho$ is $r$.

\proclaim{$\lambda$-Lemma (\cite {MSS}, \cite {BR})}
{
Let $h$ be a holomorphic motion over a hyperbolic domain
$\Lambda \subset \C$ and let $z, w \in \Lambda$.
Then $h[z,w]$ extends to a quasiconformal
map of $\C$ with dilatation bounded by $K(r)$, where $r$ is the
hyperbolic distance between $z$ and $w$ in $\Lambda$.

In the general case ($\Lambda$ not one-dimensional), the same estimate
holds with the Kobayashi distance instead of the hyperbolic distance.
In particular, if $h$ is a holomorphic motion over
$B_r(\E)$, and if $z,w \in B_{r/2}(\E)$ then $h[z,w]=1+O(\|z-w\|)$.
}

If $h=h_U$ is a holomorphic motion of an open set, we define
$\Dil(h)$ as the supremum of the dilatations of the maps
$h[z,w]$.

A {\it completion} of a holomorphic
motion means an extension of $h$ to the whole
complex plane: $X[z]=\C$ for all $z \in \Lambda$.  The problem of existence
of completions is considerably different in one-dimension or higher:

\proclaim{Extension Lemma (\cite {Sl})}
{
Any holomorphic motion over a simply connected domain $\Lambda \subset \C$
can be completed.
}

\proclaim{Canonical Extension Lemma (\cite {BR})}
{
Let $h$ be a holomorphic motion over $B_r(\E)$.  Then the restriction of $h$
to $B_{r/3}(\E)$ can be completed in a canonical way.
}

\subsubsection{Symmetry}

Let us assume that $\E$ is the complexification of a real-symmetric space
$\E^\R$, that is, there is a anti-linear isometric involution $\conj$ fixing
$\E^\R$.  Let us use $\conj$ to denote also the map $(z,w) \to (\conj
z,\overline w)$ in $\E \oplus \C$.

A set $X \subset \E, \E \oplus \C$ is called real-symmetric if $\conj(X)=X$.
Let $\Lambda \subset \E$ be a real-symmetric domain.
A holomorphic motion $h$ over $\Lambda$ 
is called real-symmetric if the image of any leaf by $\conj$ is also a leaf.

The systems we are interested on are real, so they naturally possess
symmetry.  In many cases, we will consider a real-symmetric holomorphic
motion associated to
the system, which will need to be completed using the Extension Lemma (in
one-dimension) or the Canonical Extension Lemma (in higher dimensions). 

Since the Canonical Extension Lemma is canonical, it can be used to produce
real-symmetric holomorphic motions out of real-symmetric holomorphic
motions (see Remark 2.2 of \cite {ALM}).
On the other hand, the Extension Lemma adds ambiguity on the
procedure, since the extension is
not unique.  In particular, this could lead to loss of symmetry.  In order
to avoid this problem, we will choose a little bit more carefully our
extensions.  The relevant result is then the following:


\proclaim{Real Extension Lemma.}
{
Any real-symmetric holomorphic motion over a simply connected domain
$\Lambda \subset \C$ can be completed to a real-symmetric
holomorphic motion.
}

This version of the Extension Lemma can be proved in the same way as the  
non-symmetric one\footnote {This is particularly easy to check in
Douady's proof \cite {Do} of the Extension Lemma.  Indeed, there exists
only one step which could lead to loss of symmetry, and thus needs to be
looked more carefully in order to obtain the Real Extension Lemma:
in Proposition 1 we should make sure that the
(not uniquely defined) diffeomorphism $F$ is chosen real-symmetric
(the proof that this is possible is the same).}.

So we can make the following:

\proclaim{Symmetry assumption.}
{
Extensions of real-symmetric motions will always be taken real-symmetric.
}


\subsubsection{Notation warning}

We will use the following conventions.  Instead of talking about the sets
$X[z]$, fixing some $z \in \Lambda$, we will say that $h$ is the
motion of $X$ over $\Lambda$, where $X$ is to be thought of as a set which
depends on the point $z \in \Lambda$.  In other words, we usually drop the
brackets from the notation.

We will also use the following notation for restrictions of
holomorphic motions:
if $Y \subset X$, we denote $\YY \subset \XX$ as the union of
leaves through $Y$.

\subsection{Codimension-one laminations}

Let $\F$ be a Banach space.
A codimension-one {\it holomorphic lamination} $\LL$ on an open subset
$\WW \subset \F$ is a family of disjoint codimension-one  
Banach submanifolds of $\F$, called the {\it leaves} of the lamination
such that for any point $p \in \WW$, there exists a holomorphic local chart
$\Phi:\tilde \WW \to \VV \oplus \C$,
where $\tilde \WW \subset \WW$ is a neighborhood of $p$ and
$\VV$ is an open set in some complex Banach space $\E$,
such that for any leaf $L$ and any connected component $L_0$ of $L\cap \WW$,
the image $\Phi(L_0)$  is a graph of a holomorphic function $\VV \ra \C$.

It is clear that the local theory of codimension-one laminations is the
theory of holomorphic motions.
For instance, the $\lambda$-Lemma implies that
holonomy maps of codimension-one laminations have quasiconformal extensions,
and gives bounds on the dilatation of those extensions.

\section{Complex dynamics} \label {compl dyn}

In this section we introduce some basic language necessary to describe
precisely the constructions of \cite {parapuzzle}.  Although this language
may seem at first technical and heavy, it will allow us to give formal and
concise proofs of the results we need (which are extensions of the results
of \cite {parapuzzle}).  We warn the reader that the notation is different
from \cite {parapuzzle}.


Through this section, we will deal exclusively with one-dimensional
holomorphic motions over some Jordan domain $\Lambda \subset \C$.

\comm{
The set $\Lambda$
should be thought as the parameter space of some complex family of
one-dimensional dynamical systems, and the second coordinate in $\Lambda
\oplus \C$ as the phase-space.
}

\subsection{$R$-maps and $L$-maps}

Let $U$ be a Jordan disk and $U^j$, $j \in \Z$ be a
family of Jordan disks with disjoint closures such that
$\overline {U^j} \subset U$ for every $j \in \Z$.  We assume further that
$0 \in U^0$.
A holomorphic map $R:\cup U^j \to U$ is called a $R$-map (return type map)
if for $j \neq 0$,
$R|U^j$ extends to a homeomorphism $R:\overline {U^j} \to \overline U$
and $R|U^0$ extends to a double covering map $R:\overline {U^0} \to
\overline U$ ramified at $0$.

For $\d \in \Omega$, $\d=(j_1,...,j_m)$,
we define $U^\d=\{z \in U|\, R^{k-1}(z) \in
U^{j_k},1 \leq k \leq m\}$ and we let $R^\d=R^m|U^\d$.
Let $W^\d=(R^\d)^{-1}(U^0)$.

Given an $R$-map $R$ we define an {\it $L$-map} (landing type map)
$L(R):\cup W^\d \to U^0$, by setting $L(R)|W^\d=R^\d$ (thus $L(R)$ is the
first landing map to $U^0$ under the dynamics of $R$).  We will say that
$L(R)$ is the landing map associated with (or induced from) $R$.
\comm{
Notice that
$L(R)|W^\d$ extends to a homeomorphism onto $\overline U$.
}

\subsubsection{Renormalization}

Given an $R$-map $R$ such that $R(0) \in \cup W^\d$
we can define the (generalized in the sense of Lyubich) {\it
renormalization}
$N(R)$ as the first return map to $U^0$ under the dynamics of $R$.  It
follows that $N(R)=L(R) \circ R$ where defined in $U^0$, and
that $N(R)$ is also an $R$-map.

\comm{
The $\d$ {\it truncation} of $L(R)$ is defined by
$L^\d(R)=L(R)$ outside
$U^\d$ and as $L^\d(R)=R^\d$ in $U^\d$.
Given a $R$-map $R$ such that $R(0) \in U^\d_0$
we can defined the {\it gape renormalization}
$G^\d(R)$ by $G^\d(R)=L^\d(R) \circ R$ where defined in $U^0$.
}

\comm{
Notice that either $G^\d(R)=R|U^0$ (if $\d$ is empty) or the domain of
$G^\d(R)$ consists of countably many Jordan disks (with disjoint closures)
where $G^\d(R)$ coincides with $N(R)$ and
acts as a diffeomorphism onto $U^0$ and one central Jordan disk
(containing $0$) where $G^\d(R)$ acts as a double covering onto $U$.
}

\subsection{Tubes and tube maps}

A {\it proper motion} of a set $X$ over $\Lambda$ is a
holomorphic motion of $X$ over $\Lambda$ such that for any $z \in \Lambda$,
the map $\hh[z]:\Lambda \times X[z] \to \XX$ defined by
$\hh[z](w,x)=(w,h[z,w](x))$
has an extension to $\overline \Lambda \times X[z]$
which is a homeomorphism.

An {\it equipped tube} $h_T$ is a holomorphic motion of a
Jordan curve $T$.  Its support is called a {\it tube}.
We say that an equipped tube is {\it proper} if it is a proper motion.  Its
support is called a {\it proper tube}.
The {\it filling} of a tube $\TT$ is the set
$\UU \subset \Lambda \times \C$
such that $U[z]$ is the bounded component of $\C \setminus T[z]$,
$z \in \Lambda$.

A {\it special motion} is a holomorphic motion $h=h_{X \cup T}$
such that $\XX$ is contained in the filling $\UU$ of $\TT$,
$h|T$ is an equipped proper
tube and the closure of any leaf through $X$ does not intersect the closure
of $\TT$.

If $\TT$ is a tube over $\Lambda$, and $\UU$ is its filling, a fiberwise
holomorphic map
$\FF:\UU \to \C^2$ is called a {\it tube map}
if it admits a continuous extension
to $\overline \UU$.

\subsubsection{Tube pullback}

\comm{
Let us now describe a way of defining new holomorphic motions by conformal
pullback of another holomorphic motion.
}

Let $\FF:\VV \to \C^2$ be a tube map such that
$\FF(\partial \VV)=\partial \UU$, where $\UU$ is the filling of a
tube over $\Lambda$ and
let $h$ be a holomorphic motion supported on $\overline \UU \cap
\pi_1^{-1}(\Lambda)$.
Let $\Gamma$ be a (parameter) open set such that
$\overline \Gamma \subset \Lambda$ and $W$
be a (phase) open set which moves
holomorphically by $h$ over $\Lambda$ and such that
$\overline W \subset U$.  Assume that $\overline W$ contains critical
values of $\FF|(\VV \cap \pi_1^{-1}(\overline \Gamma))$, that is,
if $\lambda \in \overline \Gamma$,
$z \in V[\lambda]$ and $DF[\lambda](z)=0$ then
$F[\lambda](z) \in \overline {W[\lambda]}$.

Let us consider a leaf of $h$ through $z \in U \setminus \overline W$,
and let us denote by $\EE(z)$ its preimage by $\FF$ intersected with
$\pi^{-1}_1(\Gamma)$.
Each connected component of $\EE(z)$ is a graph over
$\Gamma$, moreover, $\overline \EE(z) \subset \UU$.
So the set of connected components of
$\EE(z)$, $z \in U \setminus \overline W$ is a holomorphic
motion over $\Gamma$.
We define a new holomorphic motion over $\Gamma$, called
{\it the lift of $h$ by $(\FF, \Gamma, W)$}, as an extension
to the closure of $V$ of the holomorphic motion whose leaves are the
connected components of $\EE(z)$, $z \in U \setminus \overline W$
(the lift is not uniquely defined).
It is clear that this holomorphic motion
is a special motion of $V$ over $\Gamma$ and its dilatation over
$F^{-1}(U \setminus \overline W)$ is bounded by
$K(r)$ where $r$ is the hyperbolic diameter of $\Gamma$ in $\Lambda$.

\subsubsection{Diagonal and
Phase-Parameter holonomy maps} \label {holonomy maps}

Let $h$ be an equipped proper tube supported on
$\TT$.  A {\it diagonal} of $\TT$ is a holomorphic
section $\Psi:\Lambda \to \C^2$ (so that $\pi_1 \circ \Psi=\id$),
admitting a continuous extension to $\Lambda$, and such that $\Psi(\Lambda)$
is contained on the filling of $\TT$ and for $\lambda \in \Lambda$,
$h[\lambda] \circ \Psi|\partial \Lambda$ has degree one onto $T[\lambda]$.

Let $h=h_{X \cup T}$ be a special motion
and let $\Phi$ be a diagonal of $h|T$.
It is a consequence of the Argument Principle (see \cite {parapuzzle})
that the leaves of $h|X$ intersect
$\Phi(\Lambda)$ in a unique point (with multiplicity one).
From this we can define a map $\chi[\lambda]:X[\lambda] \to \Lambda$
such that
$\chi[\lambda](z)=w$ if $(\lambda,z)$ and $\Phi(w)$
belong to the same leaf of $h$.
It is clear that each $\chi[\lambda]$
is a homeomorphism onto its image, moreover,
if $U \subset X$ is open, $\chi[\lambda]|U[\lambda]$
is locally quasiconformal, and if $\Dil(h|U)<\infty$ then
$\chi[\lambda]|U[\lambda]$ is globally quasiconformal with
dilatation bounded by $\Dil(h|U)$.

We will say that
$\chi$ is the {\it holonomy family} associated to the pair $(h,\Phi)$.

\begin{rem} \label {dilatation}

Let $h_{U \cup T}$ be a special motion, $\Phi$ a diagonal,
and let $\chi$ be the holonomy family associated to $(h,\Phi)$.  Let
$X$ be compactly contained in $U$.  Then the $\lambda$-lemma implies that
for every $\lambda \in X$, $\chi[\lambda]|X[\lambda]$ extends to a qc map of
the whole plane\footnote{Here we use that the restriction of a quasiconformal
map $\chi$ to a compact subset of its domain always admits a global qc
extension (the bounds on the dilatation of the global extension depending on
the original bounds and on the hyperbolic diameter of $X$ in $U$).}.

If $\chi(X)$ has small hyperbolic diameter in
$\chi(U)$ then one can say more: this qc extension has dilatation close to
$1$.  Indeed, in this case there is a Jordan domain $X \subset U' \subset U$
such that $\chi(U')$ has small hyperbolic diameter in $\chi(U)$ and $\chi(X)$
has small hyperbolic diameter in $\chi(U')$.  Using the $\lambda$-lemma, one
sees that for $\lambda \in U'$, $\chi[\lambda]|U'[\lambda]$ has dilatation
close to $1$, and we may apply the previous argument.  (This does not work
if we only know that $X$ has small hyperbolic diameter in $U$.)

\end{rem}

\subsection{Families of $R$-maps}

An {\it $R$-family} is a pair $(\RR,h)$, where $\RR$ is a holomorphic map
$\RR:\cup \UU^j \to \UU$
such that the fibers $R[\lambda]$ of $\RR$ are
$R$-maps, for every $j$, $\RR|\UU^j$ is a tube map,
and $h=h_{\overline U}$ is a holomorphic motion
such that $h|(\partial U \cup \cup_j \partial U_j)$ is special.
If additionally $\RR \circ \0$ is a diagonal to $h$,
we say that the $\RR$ is {\it full}.

\subsubsection{From $R$-families to $L$-families}

Given an $R$-family $\RR$ with motion $h=h_{\overline U}$ we induce a family
of $L$-maps as follows.
If $\d \in \Omega$, $\d=(j_1,...,j_m)$, we let
$\UU^\d=\{(\lambda,z) \in \UU|R^{k-1}[\lambda](z) \in
U^{j_k}[\lambda]\}$ and
define $\RR^\d=\RR^m|\UU^\d$.  Let $\WW^\d=(\RR^\d)^{-1}(\UU^0)$.
We define $L(\RR):\cup \WW^\d \to \UU^0$ by
$L(\RR)|\WW^\d=\RR^\d$.  Notice that the
$L$-maps which are associated with the fibers of $\RR$ coincide with
the fibers of $L(\RR)$.

We define a holomorphic motion $L(h)$ in the following way.
The leaf through $z \in \partial U$ is the leaf of
$h$ through $z$.  If there is a smallest $U^\d$ such that
$z \in U^\d$, we let the leaf through
$z$ be the preimage by $\RR^\d$ of the leaf of $h$ through
$\RR^\d(z)$.  We finally extend it to $\overline U$ using the
Extension Lemma.

The {\it $L$-family} associated to $(\RR,h)$ is the pair $(L(\RR),L(h))$.

\subsubsection{Parameter partition and family renormalization}

Let $(\RR,h)$ be a full $R$-family.  Since $L(h)|(U \cup \cup_j \overline
{U^j})$ is special, we can consider the holonomy family of
the pair $(L(h)|(U \cup \cup_j \overline {U^j}),\RR(\0))$,
which we denote by $\chi$.  We use $\chi$ to partition $\Lambda$: we will
denote $\Lambda^\d=\chi(U^\d)$ and $\Gamma^\d=\chi(W^\d)$.


The $\d$-renormalization of
$(\RR,h)$ is the $R$-family $(N^\d(\RR),N^\d(h))$ over $\Gamma^\d$
defined as follows.  We take $N^\d(h)$ as the lift of $L(h)$ by
$(\RR|\UU^0, \Gamma^\d, W^\d)$
where defined.  It is clear that
$(N^\d(\RR),N^\d(h))$ is full,
and its fibers are renormalizations of the fibers of $(\RR,h)$.
Moreover, $N^\d(h)$ is a special motion.

\subsubsection{Chains}

An {\it $R$-chain} over $\lambda_0$
is a sequence of full $R$-families $(\RR_i,h_i)$, over
domains $\Lambda_i$, $i \geq 1$, such that $\lambda_0 \in \cap \Lambda_i$
and which are related by renormalization: $\RR_{i+1}=N^{\d_i}(\RR_i)$,
$h_{i+1}=N^{\d_i}(h_i)$ for some sequence $\d_i$.
We will say that a level $\RR_i$ of the chain is central if
$|\d_i|=0$.

\subsubsection{Gape motion}

In the situations we shall face, the central puzzle piece $U^0_i$
degenerates to a figure eight when $\lambda$ goes to the boundary of
$\partial \Lambda_i$.  This will force us to consider a technical
modification of the holomorphic motion $h_i$ as follows. 

For $i>1$, let $G(h_{i-1})$ be a holomorphic motion of $U_{i-1}$ over
$\Lambda^{\d_{i-1}}_{i-1}$ which coincides with $L(h_{i-1})$ on $U_{i-1}
\setminus \overline {U^0_{i-1}}$, and coincides with the lift of $L(h)$ by
$(\RR|\UU^0_{i-1}, \Lambda^{\d_{i-1}}_{i-1}, U^{\d_{i-1}}_{i-1})$ on
$U^0_{i-1}$.

Notice that for $i>1$, the motion $h_i$ (and hence $L(h_i)$) is special,
since it is obtained by renormalization.  So for $i>2$,
the motion $G(h_{i-1})$ is special.  Moreover, it is easy to see that
$(\RR^{|\d_{i-1}|+1}_{i-1} \circ \0)|\Lambda^{\d_{i-1}}$ (which extends
$(\RR_i \circ \0)|\Lambda_{i+1}$) is a diagonal to $G(h_{i-1})$.

\subsubsection{Real chains}

A fiberwise map $\FF:\XX \to \C^2$ is real-symmetric if $\XX$ is
real-symmetric and $\FF \circ \conj=\conj \circ \FF$.
We will say that a chain $\{\RR_i\}$ over a parameter $\lambda \in \R$
is real-symmetric if each $\RR_i$
and each underlying holomorphic motion $h_i$ is real-symmetric.

Because of the Symmetry assumption, a chain $\{\RR_i\}$ over a
parameter $\lambda \in \R$ is real-symmetric
provided the first step data $\RR_1$ and $h_1$ is real-symmetric.
In this case, all objects related to the chain are real-symmetric.

\begin{rem}

If $\RR_1$ is real-symmetric then $h_1$ can always be modified
to be real-symmetric.  Indeed if $\RR_1$ is real-symmetric then
$\partial \UU_1 \cup \cup \partial \UU^j_1$ is a real-symmetric set
and it is enough to check that $h_1|(\partial U_1 \cup \cup \partial U^j_1)$
is already real-symmetric.  To see this, first notice that
if $X[\lambda]$ moves holomorphically and $X$ has empty interior then the
holomorphic motion of $X$ is unique\footnote {In this case $\C \setminus X$
also moves holomorphically by some
motion $h_{\C \setminus X}$ obtained from the Extension Lemma, and the
motion of $X$ can be seen as coming from the extension of $h_{\C \setminus
X}$ to the closure $\overline {\C \setminus X}=\C$, and this extension is
unique.}.  This implies that if $X[\lambda]$ is real symmetric with empty
interior then any motion of $X[\lambda]$ is also real-symmetric.

\end{rem}

\subsection{Complex Phase-Parameter estimates}

We shall now show how estimates on the geometry of parapuzzle pieces yield
automatically estimates on the regularity of holonomy maps.  We shall need
four specific statements, contained in two lemmas.

\begin{lemma} \label {cphpa1}

Let us consider an $R$-chain $(\RR_i,h_i)$ over $\lambda_0$, and
let $\tau_i$ be such $R_i[\lambda_0](0) \in U^{\tau_i}_i[\lambda_0]$.
For $i>1$, let $\chi_i$ be the holonomy family associated to
$(L(h_i),\RR_i \circ \0)$.
For every $\g>1$ there exists $K>0$ such that if $\mod(\Lambda_{i-1}
\setminus \overline {\Lambda_i})>K$ and $\mod(U_{i-1}[\lambda_0] \setminus
\overline {U_i[\lambda_0]})>K$ then

\noindent {\rm [CPhPh1]} For every $\lambda \in \Lambda^{\tau_i}_i[\lambda_0]$,
$L(h_i)[\lambda_0,\lambda]|U_i[\lambda_0]$ has a $\g$-quasiconformal
extension to the whole complex plane,

\noindent {\rm [CPhPa1]}
$\chi_i[\lambda_0]|U^{\tau_i}_i[\lambda_0]$ has a $\g$-quasiconformal
extension to the whole complex plane.

Moreover, if the $R$-chain
$(\RR_i,h_i)$ is real then the claimed extensions can be
taken real as well.

\end{lemma}

\begin{pf}

Both items (1) and (2) follow easily from the $\lambda$-lemma (see also
Remark~\ref {dilatation}) if we can
stablish that $\mod(\Lambda_i \setminus \overline {\Lambda^{\tau_i}_i})$
is big.

The hypothesis on $\mod(U_{i-1}[\lambda_0] \setminus \overline
{U_i[\lambda_0]})$ implies that
$\mod(U_i[\lambda_0] \setminus \overline {U^{\tau_i}_i[\lambda_0]})$ and
$\mod(U_i[\lambda_0] \setminus \overline {U^0_i[\lambda_0]})$ are bigger
than $K/2$.  If $K$ is big, this implies that there is an annulus of big
modulus contained in $U_i[\lambda_0] \setminus \overline {U^0_i[\lambda_0]}$
and going around $U^{\tau_i}_i[\lambda_0]$.  Using again that $K$ is big and
the hypothesis on $\mod(\Lambda_{i-1}
\setminus \overline {\Lambda_i})$, we
see that the dilatation of $h_i|(U_i \setminus \overline {U^0_i})$
is small ($\lambda$-lemma).  We conclude that $\mod(\Lambda_i \setminus
\overline {\Lambda^{\tau_i}_i})$ is big as required.
\end{pf}

\begin{lemma} \label {cphpa2}

Let us consider an $R$-chain $(\RR_i,h_i)$ over $\lambda_0$.  For $i>2$,
let $\tilde \chi_i$ be the holonomy family associated to
$(G(h_{i-1}),\RR^{|\d_{i-1}|+1}_{i-1} \circ \0)$
For every $\g>1$ there exists $K>0$ such that if $\mod(\Lambda_{i-2}
\setminus \overline {\Lambda_{i-1}})>K$ and $\mod(U_{i-1}[\lambda_0]
\setminus \overline {U_{i-1}[\lambda_0]})>K$ then

\noindent {\rm [CPhPh2]}
For every $\lambda \in \Lambda_i$,
$G(h_{i-1})[\lambda_0,\lambda]|U_{i-1}[\lambda_0]$ has a $\g$-quasiconformal
extension to the whole complex plane,

\noindent {\rm [CPhPa2]}
$\tilde \chi_i[\lambda_0]|U_i[\lambda_0]$ has a $\g$-quasiconformal
extension to the whole complex plane.

Moreover, if the $R$-chain $(\RR_i,h_i)$
is real then the claimed extensions can be taken real as well.

\end{lemma}

\begin{pf}

Both items (1) and (2) follow easily from the $\lambda$-lemma (see also
Remark~\ref {dilatation}) if we can
stablish that $\mod(\Lambda^{\d_{i-1}}_{i-1} \setminus \overline
{\Lambda_i})$ is big.

The hypothesis on $\mod(\Lambda_{i-2} \setminus \overline {\Lambda_{i-1}}$
implies that the dilatation of $L(h_{i-1})|(U_{i-1} \setminus \overline
{U_i})$ is less than $2$ (provided $K$ is sufficiently big).  Notice that
$\Lambda^{\d_{i-1}}_{i-1} \setminus \overline
{\Lambda_i}=\chi_{i-1}[\lambda_0](U^{\d_{i-1}}_{i-1}[\lambda_0]
\setminus \overline {W^{\d_{i-1}}_{i-1}[\lambda_0]})$,
where $\chi_{i-1}$ is the holonomy family
associated to $(L(h_{i-1}),\RR_{i-1} \circ \0)$.
The hypothesis on $\mod(U_{i-1}[\lambda_0] \setminus \overline
{U^0_{i-1}[\lambda_0]})$ (which equals $\mod(U^{\d_{i-1}}_{i-1}[\lambda_0]
\setminus \overline {W^{\d_{i-1}}_{i-1}[\lambda_0]})$) then implies that
$\mod(\Lambda^{\d_{i-1}}_{i-1} \setminus \overline {\Lambda_i})$ is big (at
least $K/2$) as required.
\end{pf}

\comm{
Before stating the next lemma, we shall need to define an additional 

\begin{lemma}

andsince
$\mod(U_{i-1}[\lambda_0] \setminus \overline {U^0_{i-1}[\lambda_0]})$ is
big, both $\mod(U_i[\lambda_0

}

\section{Puzzle and parapuzzle geometry} \label {puzzle}

In this section we will recall an important part of Lyubich's theory of the
quadratic family (regarding linear growth of moduli of certain phase and
parameter annuli), and will discuss the validity of those results in the
context of more general $R$-chains.

\subsection{Puzzle estimates}

The following result is contained on (the proof of)
Theorem II of \cite {puzzle}:

\begin{thm} \label {3.3}

For every $C>0$, there exists $C'>0$
with the following property.  Let $R_i$ be a sequence of $R$-maps such that
$R_{i+1}=N(R_i)$ and let $n_k-1$ be the sequence of non-central
levels, so that $R_{n_k-1}(0) \notin U^0_{n_k-1}$.
If $\mod(U_1 \setminus \overline {U^0_1})>C'$ then
$\mod (U_{n_k} \setminus \overline {U^0_{n_k}})>C$.

\end{thm}

(In Lyubich's notation, $R$-maps are called generalized quadratic maps.)


The following result is Theorem III of \cite {puzzle}:

\begin{thm} \label {3.4}

For every $C'>0$, there exists $C''>0$
with the following property.  Let $R_i$ be a sequence of $R$-maps such that
$R_{i+1}=N(R_i)$ and let $n_k-1$ be the sequence of non-central
levels.  If $\mod(U_1 \setminus \overline {U^0_1})>C'$ then
$\mod (U_{n_k} \setminus \overline {U^0_{n_k}})>C''k$.

\end{thm}

\comm{
Of course both theorems imply the following:

\begin{thm}

There exists a constant $C''>0$ such that for every $C'>0$, there exists
$C>0$, with the following property.  Let $R_i$ be a sequence of
$R$-maps such that $R_{i+1}=N(R_i)$ and let $n_k-1$ be the sequence
of non-central levels.  If $\mod(U_1 \setminus \overline {U^0_1})>C'$ then
$\mod (U_{n_k} \setminus \overline {U^0_{n_k}})>\max \{C''k-C^{-1}, C\}$.

\end{thm}
}

\subsection{Parapuzzle estimates}

\subsubsection{The quadratic family}

Let $p_c(z)=z^2+c$ be the quadratic family.  The following result is
contained in Lemma 3.6 of \cite {parapuzzle}:

\begin{thm} \label {p_lambda}

Let us fix a non-renormalizable quadratic polynomial $p_{c_0}$
with a recurrent critical point and no neutral periodic orbits.
Then there exists a full
$R$-family $\RR_1$ over some $c_0 \in \Lambda_1$ such that if $c \in
\Lambda_1$ then $R[c]:\cup U^j_1[c] \to U_1[c]$ is the first return map
under iteration by $p_c$.

\end{thm}

The following is Theorem A of \cite {parapuzzle}:

\begin{thm} \label {3.5}

In the setting of Theorem \ref {p_lambda}, let $\RR_i$ be the $R$-chain over
$c_0$ with first step $\RR_1$.  If $n_k-1$ denotes the $k$-th non-central
return, then $\mod (\Lambda_{n_k} \setminus \overline {\Lambda_{n_k+1}}) >
T k$, for some
constant $T>0$.

\end{thm}

\begin{rem}

In Lyubich's notation he lets $\Delta^i=\Lambda_{n_i}$ and
$\Pi^i=\Lambda^0_{n_i}$.  He states that both $\mod (\Delta^i \setminus
\overline {\Delta^{i+1}})$ and $\mod (\Delta^i \setminus \overline {\Pi^i})$
grow linearly.  His statement implies ours after one notices that
if $n_i+1=n_{i+1}$ then $\Delta^{i+1}=\Lambda_{n_i+1}$,
otherwise $\Pi^i=\Lambda_{n_i+1}$.

\end{rem}

Those two results are proved in a slightly more general setting then we
state here: they are valid for so-called full unfolded families of
quadratic-like maps.  This version allows one to state results also for
finitely renormalizable quadratic polynomials (via renormalization).

\subsubsection{General case}

The following more general theorem can be proved using
the ideas of Theorem A of
\cite {parapuzzle} but it is a little bit tedious to check the details
(it is necessary to get deep into the construction of \cite {puzzle}).

\begin{thm} \label {fullthm}

For every $K>1$, $T>0$, there exists $T'>0$
with the following property.  Let $(\RR_i,h_i)$ be a $R$-chain over
$\lambda_0$ and let $n_k-1$ be the sequence of non-central levels.
If $\Dil(h_1|(U_1 \setminus \overline {U^0_1}))<K$ and
$\mod(U_1[\lambda] \setminus \overline {U^0_1[\lambda]})>T$ then
$\mod (\Lambda_{n_k} \setminus \overline {\Lambda_{n_k+1}})>T' k$.

\end{thm}

Since we do not need the full strength of the previous theorem, we will
state and prove a weaker estimate using a simple inductive argument.


\comm{
\begin{lemma} \label {3.7}

There exists a constant $C_0$ with the following property.
Let $h$ be a holomorphic motion over $\Lambda$ and let $\Lambda' \subset
\Lambda$ be such that $\nu=\mod (\Lambda \setminus \overline {\Lambda'}) >
C_0$.  Then the dilatation of $h$ over $\Lambda'$ is bounded by $2$.

\end{lemma}
}

\begin{thm} \label {3.8}

For every $K>1$, there exists constants $T'>0$, $T''>0$
with the following properties.  Let $(\RR_i,h_i)$ be a $R$-chain over
$\lambda_0$ and let $n_k-1$ be the sequence of non-central levels.
If $\Dil(h_1|(U_1 \setminus \overline {U^0_1})<K$ and
$\mod(U_1[\lambda_0] \setminus \overline {U^0_1[\lambda_0]})>T'$ then
$\mod (\Lambda_{n_k} \setminus \overline {\Lambda_{n_k+1}})>T'' k$.

\end{thm}

\begin{pf}

Let $\nu_i=\mod(\Lambda_i \setminus \overline {\Lambda_{i+1}})$,
$\mu_i=\mod(U_i[\lambda_0] \setminus \overline {U_{i+1}[\lambda_0]})$,
$k_i=\Dil(h_i|U_i[\lambda_0] \setminus \overline {U_{i+1}[\lambda_0]})$.
For $i>1$, denote by $\chi^0_i$ the holonomy family associated to
$(h_i,\RR_i \circ \0)$.

\comm{
We have the
following simple rules: $\mu_{i+1} \geq \mu_i/2$, $\nu_i \geq \mu_i/k_i$,
}
Notice that if $\nu_i>C_0=1000$ then $k_{i+1} \leq 2$.
Moreover, for $i>1$, and in particular for
$i=n_k$, we have $\nu_i=\mod(\chi^0_i(U_i[\lambda_0] \setminus \overline
{U_{i+1}[\lambda_0]})) \geq \mu_i/k_i$.

Let $T=(4+2C_0) (2+K)$ and let $T'>T$ be so big that if $\mu_1>T'$ then
$\mu_{n_k}>T$, $k \geq 1$.  Let also $T''$ be such that
if $\mu_1>T'$ then $\mu_{n_k}>k T''(2+K)$.

Let us assume that for some $m$, we have $\mu_m>T$ and $k_m \leq (2+K)$,
and let $m' \geq m$ be the next non-central return.

For $\lambda \in \Lambda_{m'+1}$, we have
$R_m^{m'-m+1}(0) \in W^\d_m$ for some $\d$.  Let $\Upsilon$ be
the component of $R_m(0)$ of $(R_m^{m'-m}|U_{m'})^{-1}(U^\d_m)$ and
$\Upsilon'$ be the component of $R_m(0)$ of
$(R_m^{m'-m}|U_{m'})^{-1}(W^\d_m)$.

Let $H_m=L(h_m)$ outside of $U^\d_m$ and let the leaves of
$H_m|U^\d_m$ be the preimages by $R^\d_m$ of the leaves of $h_m$.
If $m=m'$, let $H=H_m$.  Otherwise, notice that if
$\lambda \in \Lambda_{m'-1}$, then $R_m^{m'-m}|U^0_{m'-1}[\lambda]$
is a $2^{m'-m}$ branched
covering map over $U_m[\lambda]$, and for $\lambda \in \Lambda_{m'}$,
$R_m^{m'-m}|U_{m'}[\lambda] \setminus
\overline {U^0_{m'}[\lambda]}$
is unbranched.  Let $H$ be the lift of
$H_m$ by $(\RR_m^{m'-m}|U_{m'},U^0_{m'},\Lambda_{m'})$.  So in both
cases, $H$ is a holomorphic motion over $\Lambda_{m'}$.

With this definition, $\Upsilon$ and $\Upsilon'$
(which are apriori defined over $\Lambda_{m'+1}$)
move holomorphically with $H$ (over $\Lambda_{m'}$).

Let $\chi$ be the holonomy family of the pair
$(\RR_{m'} \circ \0,H|\partial U_{m'} \cup \overline \Upsilon)$.
It is clear that $\Dil(\chi|\Upsilon)$ is bounded by $k_m$.
In particular, we can estimate $\nu_{m'} \geq \mod (\Upsilon[\lambda_0]
\setminus \overline {\Upsilon'[\lambda_0]})/k_m=\mu_m/k_m \geq
\mu_m/(2+K)>C_0$.
With $m=1$, we have $k_1 \leq K \leq 2+K$ by hypothesis and
$m'=n_1-1$, so $\nu_{n_1-1} \geq \mu_1/(2+K) \geq T/(2+K) \geq C_0$ and
$k_{n_1} \leq 2 \leq 2+K$.
With $m=n_k$, we have that $m'=n_{k+1}-1$ and
$\nu_{n_{k+1}-1} \geq \mu_{n_k}/(2+K) \geq T/(2+K) \geq C_0$ and
$k_{n_{k+1}} \leq 2 \leq 2+K$, provided $k_{n_k} \leq 2+K$.
By induction, we have $k_{n_k} \leq 2+K$ for every $k$, so $\nu_{n_k} \geq
\mu_{n_k}/(2+K)>T''k$.
\end{pf}

This simple inductive argument does not seem to work easily to get the full
Theorem \ref {fullthm}\footnote {However, we will see that this
is enough to yield the full power of Theorem B of \cite
{parapuzzle} (almost every non-regular
finitely renormalizable quadratic map is
stochastic), through the arguments of this paper.  This approach only uses
geometric estimates of puzzle pieces for {\it real maps}, and may be
useful for generalizations beyond unimodal maps with a quadratic critical
point.}.

\section{Unimodal maps} \label {basic}

We refer to the book of de Melo \& van Strien \cite{MS}
for the general background in one-dimensional dynamics.

We will say that a smooth (at least $C^2$) map $f: I\ra I$ of the
interval $I=[-1,1]$ is {\it unimodal} if $f(-1)=-1$, $f(x)=f(-x)$ and $0$ is
the only critical point of $f$ and is non-degenerate, so that
$D^2 f(0) \neq 0$.

\begin{rem}

The introduction of normalization and symmetry assumptions was made
in order to avoid cumbersome notations: all results and proofs
generalize to the non-symmetric case.  See also Appendix C of \cite {ALM}.

\end{rem}

\begin{rem}

The assumption that the critical point is non-degenerate is made just for
convenience: typical unimodal maps certainly have non-degenerate
critical point.  If one is not willing to make this assumption already in
the definition, one should
add the non-degeneracy condition to the Kupka-Smale definition below.
In this case it would still hold that in non-trivial analytic families
the set of
parameters with a degenerate critical point have zero Lebesgue measure (and
is contained in a countable number of analytic subvarieties with
codimension at least $1$), see Lemma \ref {ksmany}.

The theory of unimodal maps with fixed non-quadratic criticality is
considerably different and less complete than the typical case,
and the proofs of this work do not apply.

\end{rem}

Let $\U^k$, $k \geq 2$ be the space of $C^k$ unimodal maps.
We endow $\U^k$ with the $C^k$ topology.

Basic examples of unimodal maps are given by quadratic maps
\begin{equation}\label{quadratic family}
    q_\tau \colon I \to I,\quad  q_\tau(x)= \tau-1-\tau x^2,  
\end{equation}
where $\tau\in [1/2,2]$ is a real parameter.

A map $f \in \U^2$ is said to be Kupka-Smale if all periodic orbits are
hyperbolic.  It is said to be hyperbolic if it is Kupka-Smale and
the critical point is attracted to a periodic attractor.  It is said to be
regular if it is hyperbolic and its critical point is not periodic or
preperiodic.  It is well known that regular maps are structurally stable.

A $k$-parameter $C^r$ (or analytic) family of unimodal maps
is a $C^r$ (or analytic) map $F:\overline \Lambda \times I \to I$
such that $f_\lambda \in \U^2$, where $f_\lambda(x)=F(\lambda,x)$ where
$\Lambda \subset \R^k$ is a bounded open connected domain
with smooth ($C^\infty$) boundary.  We denote $\UF^r(\Lambda)$ the space of
$C^r$ families of unimodal maps, endowed with the $C^r$ topology.
Notice that $\UF^r(\Lambda)$ is a separable Baire space.

We will not introduce a topology in the space of
analytic families of unimodal maps.

\subsection{Combinatorics and hyperbolicity}

Let $f \in \U^2$.  A symmetric interval $T \subset I$
is said to be nice if the iterates of $\partial T$ never return to
$\inter T$.  A nice interval $T \neq I$ is said to be a
restrictive (or periodic) interval of period $m$ for $f$ if
$f^m(T) \subset T$ and $m$ is
minimal with this property.  In this case,
the map $A \circ f^m \circ A^{-1}:I \to I$ is again unimodal for some affine
map $A:T \to I$: this map is usually called a
renormalization of $f$ if $m>1$ or a unimodal restriction if $m=1$.

If $T \subset I$ is a nice interval, the domain of the
first return map $R_T$ to $T$ consists of a (at most)
countable union of intervals which we denote $T^j$.  We reserve the
index $0$ for the component of $0$: $0 \in T^0$, if $0$ returns to $T$.
From the nice condition, $R_T|T^j$ is a diffeomorphism if
$0 \notin T^j$, and is an even map if $0 \in T^j$.  We call
$T^0$ the central domain of $R_T$.  The return $R_T$ is said
to be central if $R_T(0) \in T^0$.

Under the Kupka-Smale condition, the dynamics outside a nice interval
is hyperbolic, and in particular persistent:

\begin{lemma} \label {5.2}

Let $f \in \U^2$ and let $T \subset I$ be a symmetric interval.
If all periodic orbits contained if $I \setminus \inter T$ are
hyperbolic (in particular if $f$ is Kupka-Smale), then

(1)\, The set of points $X \subset I$
which never enter $\inter T$ splits in two forward invariant sets:
an open set $U$ attracted by a finite number of periodic orbits
and a closed set $K$ such that $f|K$ is uniformly expanding:
$|Df^n(x)|>C \lambda^n$, for $x \in K$ and for some constants $C>0$,
$\lambda>1$.  Moreover, preperiodic points are dense in $K$.

(2)\, There exists a neighborhood $\VV \subset \U^2$ of $f$ and a
continuous family of homeomorphisms $H[g]:I \to I$, $g \in \VV$
such that $g \circ H[g]|I \setminus T=H[g] \circ f$, and $H[f]=\id$.

\end{lemma}

\begin{pf}

The first item is a consequence of Ma\~n\'e's Theorem (see \cite {MS},
Theorem 5.1 and Corollary 1, page 248).  Since hyperbolic sets are
persistent, the second item follows.
\end{pf}


The following well known result shows that nice intervals allow one
to study arbitrarily small neighborhoods of $0$.

\begin{lemma} \label {5.1}

Let $f \in \U^2$ be Kupka-Smale.  If $f$ is not hyperbolic and the critical
orbit is infinite, then for every $\epsilon>0$,
there exists a nice interval $[-p,p] \subset (-\epsilon,\epsilon)$
with $p$ preperiodic.

\end{lemma}

\begin{pf}

Let $T$ be the intersection of all nice intervals containing $0$ whose
boundary is preperiodic.  If $T \neq \{0\}$, then the domain of $R_T$ is
either $T$ or empty.  In the first case, $R_T:T \to T$ has no fixed point in
$\inter T$ and it follows that $R_T^m(\inter T)$ converge to a periodic
attractor in $\partial T$.  Otherwise, by Lemma \ref {5.2},
$\inter f(T)$ must be contained in the basin of a periodic attractor, so
$f$ is either hyperbolic or the critical point is preperiodic.
\end{pf}

\comm{
\begin{pf}

If $0$ is recurrent, one considers either renormalization or the principal
nest.  Otherwise, $f(0)$ lands in a hyperbolic set.
\end{pf}
}


The following is an easy consequence of Lemma \ref {5.2}.

\begin{lemma} \label {5.3}

Let $f_\lambda$, $\lambda \in (-\epsilon,\epsilon)$
be a $C^2$ family of unimodal maps, and let $T$ be
a nice interval with preperiodic boundary for $f=f_0$.
Assume that there exists an interval $0 \in J$
and a family $T[\lambda]$ of intervals with
preperiodic boundary,
such that $T[0]=T$ and for $\lambda \in J$, all non-hyperbolic periodic orbits
of $f_\lambda$ intersect
$\inter T[\lambda]$.  Then there exists a continuous family of homeomorphisms
$H[\lambda]:I \to I$, $\lambda \in J$
such that $H[\lambda](T)=T[\lambda]$ and
$f_\lambda \circ H[\lambda]|(I \setminus T)=H[\lambda] \circ f$ and
$H[0]=\id$.

\end{lemma}



\subsubsection{Principal nest}

We say that $f$ is infinitely renormalizable if there exists arbitrarily
small restrictive intervals $T \subset I$.  Otherwise we say that $f$ is
finitely renormalizable.

Let $\FF \subset \U^2$ be the class of
Kupka-Smale finitely renormalizable maps whose critical point is
recurrent, but not periodic.  If $f \in \FF$, the first return map
$f^m:T \to T$ to its smallest restrictive interval has an orientation
reversing fixed point which we call $p$.  Let $I_1=[-p,p]$.  Define a nested
sequence of intervals $I_i$ as follows.  Assuming $I_i$ defined,
let $R_i$ be the first return map to $I_i$ and let $I_{i+1}$ be the central
domain $I^0_i$ of $R_i$.

The sequence $I_i$ is called the {\it principal nest} of $f$.  A level of
the principal nest is called central if $R_i$ is a central return.  We say
that a map $f \in \FF$ is simple if there are only finitely many non-central
levels in the principal nest.

\comm{
Applying Lemma \ref {5.2} we get:

\begin{lemma}

Let $f \in \FF$.  For every $i>0$, there exists a neighborhood
$\tilde \VV_i \subset \U^2$ of $f$ such that for every $g \in \VV_i$
there exists an interval $I_i[g]$ and a map $H[g]:I \to I$ such that
$g \circ H[g]|(I \setminus I_i)=H[g] \circ f$.  The map $H[g]$ depends
continuously on $g$.

\end{lemma}
}

\subsection{Negative Schwarzian derivative}

The Schwarzian derivative of a $C^3$ map $f:I \to I$ is defined by
$$
Sf=\frac {D^3 f} {Df}-\frac {3} {2} \left (\frac {D^2 f} {Df} \right )^2
$$
in the complement of the critical points of $f$.  If $Sf$ and $Sg$ are
simultaneously positive (or negative) then
$S(g \circ f)$ is positive (or negative).

If $f$ is a unimodal map the condition of negative Schwarzian derivative is
very useful and can be exploited in several ways.  One of the most used
tools is the Koebe Principle:

\begin{lemma}[Koebe Principle, see \cite {MS}, page 258]

Let $f:T \to \R$ be a diffeomorphism with non-negative Schwarzian derivative. 
Then for every $K_0$, there exists a constant $k_0$ such that if
$T' \subset T$ and both components $L$ and $R$ of $T \setminus T'$ are
bigger than $K |T'|$ for some constant $K>K_0$
then the distortion of $f|T'$ is bounded by $k_0$.
In particular, we have
$\min \{|f(L)|,|f(R)|\} \geq \hat k_0 K |f(T')|$, for some $\hat
k_0$ depending only on $K_0$.  Moreover, $k_0 \to 1$ as $K_0 \to \infty$.

\end{lemma}

Quadratic maps have
negative Schwarzian derivative.  Moreover, one can often reduce to this
situation as is shown by the following well known estimate:

\begin{lemma} \label {5.6}

If $f \in \U^3$ is infinitely renormalizable, then if $T \subset I$ is a
small enough periodic nice interval, the first return map to $T$ has
negative Schwarzian derivative.

\end{lemma}

Recently, Kozlovski showed that the assumption of negative Schwarzian can be
often removed.  The next result follows from Lemma \ref {5.2}
and \cite {GSS} (which is based on the work of Kozlovski \cite {K}).

\begin{lemma} \label {5.7}

Let $f \in \FF \cap \U^3$.  There exists $i>0$, an analytic
diffeomorphism $s:I \to I$ and a neighborhood
$\VV \subset \U^3$ of $f$, such that
there exists a continuation $I_i[g]$, $g \in \VV$ of $I_i$
($H[g](I_i)=I_i[g]$ in the notation of Lemma \ref {5.2}) such that
the first return map to $s(I_i[g])$ by $s \circ g \circ s^{-1}:I \to I$
has negative Schwarzian derivative.

\end{lemma}

\subsection{Decay of geometry}

The following result is due to Lyubich \cite {attractors}
in the case of negative Schwarzian
derivative and holds in general due to the work of Kozlovski:

\begin{lemma} \label {5.5}

Let $f \in \FF$ be at least $C^3$, and let $n_k-1$ denote the sequence of
non-central levels in the principal nest of $f$.
Then $|I_{n_k+1}|/|I_{n_k}|<C\lambda^k$ for
some constants $C>0$, $\lambda<1$.

\end{lemma}

\subsection{Quasiquadratic maps}

A map $f \in \U^3$ is {\it quasiquadratic} if any nearby map $g \in
\U^3$ is topologically conjugate to some quadratic map.
By the theory of Milnor-Thurston and Guckenheimer \cite{MS},
a map $f\in \U^3$ with negative  Schwarzian derivative and
$D^2 f(-1)<0$ is quasiquadratic, so quadratic maps are quasiquadratic.
The following results give sufficient conditions for a unimodal map to be
quasiquadratic:

\begin{thm}[see Lemma 2.13 of \cite {ALM}] \label {5.8}

Let $f \in \U^3$ be a Kupka-Smale unimodal map which is topologically
conjugate to a quadratic map.  Then $f$ is quasiquadratic.

\end{thm}

\begin{thm}[see Remark 2.6 of \cite {ALM}] \label {quad}

Let $f \in \U^3$.  If $f$ is not conjugate to a quadratic polynomial then
there exists a (not necessarily hyperbolic)
periodic orbit which attracts an open set.
In particular, if all periodic orbits of $f$ are repelling then $f$ is
quasiquadratic.

\end{thm}

\begin{rem}

Theorem \ref {5.8} is the reason that the quasiquadratic condition
considers only $C^3$ maps and the $C^3$ topology (otherwise it would not be
possible to guarantee that even quadratic maps are quasiquadratic).

\end{rem}

\subsection{Spaces of analytic unimodal maps}

Let $a>0$, and let $\Omega_a \subset \C$ be the set of points at distance at
most $a$ of $I$.  Let $\EE_a$ be the
complex Banach space of holomorphic maps
$v:\Om_a \to \C$ continuous up to the boundary
which are $0$-symmetric (that is, $v(z)=v(-z)$)
and such that $v(-1)=v(1)=0$,
endowed with the $\sup$-norm $\|v\|_a=\|v\|_\infty$.
It contains the real Banach space $\EE_a^{\R}$ of ``real maps'' $v$, 
i.e, holomorphic maps symmetric with respect to the real line: 
$v(\overline z) = \overline {v(z)}$. 

Let us consider the constant function $-1 \in \Omega_a$.
The  complex affine subspace $-1+\EE_a$ will be denoted as $\AAA_a$.

Let $\U_a=\U^2 \cap \AAA_a$.
It is clear that any analytic unimodal map belongs to some $\U_a$.
Note that $\U_a$ is the union of an
open set in the affine subspace $\AAA_a^\R= -1+ \EE_a^\R$ and a
codimension-one space of unimodal maps satisfying $f(0)=1$.

\subsection{Hybrid lamination} \label {hybrid}

One of the main results of \cite {ALM} is to describe the structure of
the partition in topological classes of spaces of analytic unimodal maps.
In that paper, they consider only the quasiquadratic case, but their proof
works for the general case (due to the results of Kozlovski) and gives the
following:

\begin{thm}[Theorem A of \cite {ALM}] \label {8.1}

Let $f \in \U_a$ be a Kupka-Smale map.  There exists a
neighborhood $\VV \subset \AAA_a$ of $f$ endowed with a codimension-one
holomorphic lamination $\LL$ (also called hybrid lamination)
with the following properties:

(1)\, the lamination is real-symmetric;

(2)\, if $g\in \VV\cap \AAA_a^\R$ is non-regular,
then the intersection of the leaf through 
$g$ with $\AAA_a^\R$ coincides with the intersection of the topological
conjugacy class of $g$ with $\VV$;

(3)\, Each $g \in \VV \cap \AAA^\R_a$ belongs to some leaf of $\LL$.

\end{thm}

(For the definition of the leaves of $\LL$ in the regular case, see
Appendix A.)

\begin{thm} \label {8.2'}

In the setting of Theorem \ref {8.1}, if $g_1,g_2 \in \VV$ are in the same
leaf of $\LL$ and $\g_1(\lambda)$, $\g_2(\lambda)$ are real
analytic paths in $\VV \cap \AAA^\R_a$, transverse to the leaves of $\VV$
and such that $\g_1(\lambda_1)=g_1$, $\g_2(\lambda_2)=g_2$, then the
local holonomy map
$\psi:(\lambda_1-\epsilon,\lambda_1+\epsilon) \to
(\lambda_2-\epsilon',\lambda_2+\epsilon')$ is quasisymmetric.  Moreover,
for $\delta$ sufficiently small,
$\psi|(\lambda_1-\delta,\lambda_1+\delta)$ is $1+O(\|g_1-g_2\|_a)$-qs.

\end{thm}

\begin{pf}

This estimate
is just the $\lambda$-Lemma in the context of codimension-one
complex laminations.
\end{pf}


Moreover, each non-regular topological class is like a Teichmuller space:

\begin{thm} \label {8.2}

In the setting of Theorem \ref {8.1},
if $g_1,g_2 \in \VV \cap \U_a$ belong to the
same leaf of $\LL$, then there exists a $1+O(\|g_1-g_2\|_a)$-qs map
$h:I \to I$ such that $g_2 \circ h=h \circ g_1$.

\end{thm}

\begin{pf}

This follows from Proposition 8.9 of \cite {ALM} and the $\lambda$-Lemma.
\end{pf}

The tangent space to topological classes has a nice characterization:

\begin{thm}[Theorem 8.10 of \cite {ALM}] \label {8.4}

If $f \in \U_a$ is a non-regular Kupka-Smale map then the
tangent space to the topological class of $f$ is given by the set of vector
fields $v \in \EE_a$ which do not admit a representation $v=\alpha \circ
f-\alpha Df$ on the critical orbit with $\alpha$ a qc vector field of $\C$.

\end{thm}

\subsection{Analytic families}

Let $\{f_\lambda\}_{\lambda \in \Lambda}$ be an analytic family of unimodal
maps.  Then for $a>0$ sufficiently small, $\lambda \mapsto f_\lambda$ is an
analytic map from $\Lambda$ to $\U_a$.  We say that $f_\lambda$ is
{\it non-trivial} if the set of regular parameters is dense.

If $\lambda_0 \in \Lambda$ is a
Kupka-Smale parameter, transversality to the topological class of
$\lambda_0$ has the obvious meaning (using Theorem \ref {8.1}).  We remark
that this definition does not depend on the choice of $\U_a$.

\begin{rem} \label {dense}

Let $B_i$ be an enumeration of all open balls contained in $\Lambda$ of
rational radius and center.  The condition of non-triviality of a family
$\{f_\lambda\}$, $\lambda \in \Lambda$ is an intersection of a countable
number of conditions (existence of a regular parameter $\lambda \in B_i$). 
Each of those conditions is open in $\UF^2(\Lambda)$.
The set of non-trivial analytic
families is also dense in the $\UF^\infty(\Lambda)$ (this would still hold
natural topology of analytic families in $\Lambda$,
which we did not introduce), due to Theorem \ref {8.1}.

We should remark that for an analytic family of quasiquadratic maps,
non-triviality is equivalent to existence of {\it one} regular parameter
(since all non-regular topological classes are analytic submanifolds in the
quasiquadratic case).  In particular, non-triviality is a $C^3$ open
condition in the quasiquadratic case.

\end{rem}

\section{Construction of the special family} \label {family}

\subsection{Puzzle maps} \label {6.1}

Let $f \in \U_a$ be a finitely renormalizable unimodal map with a
recurrent critical point.  Let us consider some nice interval
$A^0$ and let $\{A^j\}$ be the connected components of the domain of the
first landing map from $I$ to $A^0$.  We call the family
$\{A^j\}$ the real puzzle for $f$ associated to $A^0$.  The basic
object used in \cite {ALM} to analyze the dynamics
of unimodal maps can be viewed as a complexification of such
real puzzles, which are called simply a puzzle.

The definition of puzzle in \cite {ALM} is too general and technical for our
purposes.  In this paper, we will simply describe how to construct
a puzzle for $f$ (or rather a geometric puzzle, in the language of
\cite {ALM}).  Instead of giving the precise definitions of a puzzle, we
will just obtain the properties that are needed for our results.

Let us fix some advanced level $\n$ of the principal nest
of $f$ and assume that $|I_\n|/|I_{\n-1}|$ is very small.
Let us fix the following notation: let $A^0=I_\n$ and let $\{A^j\}$
be the real puzzle associated to $A^0$.  We let $A^1$ be such that
$f(0) \in A^1$.

Given $0<\theta \leq \pi/2$, and $A \subset \R$, let $D_\theta(A)$ be the
intersection of two round disks $D_1$ and $D_2$ where $D_1 \cap
\R=A$, $\partial D_1$ intersects $\R$ making an angle $\theta$, and $D_2$ is
the image of $D_1$ by symmetry about $\R$.  The
complexification of the real puzzle $\{A^j\}$ should be imagined as
$\{D_\theta(A^j)\}$ for a suitable value of $\theta$.  Of course, since the
system is non-linear, the definition can not be so simple.
Nevertheless, the
condition $|I_\n|/|I_{\n-1}|$ small allows one to bound the
nonlinearity of the
first landing map to $I_\n$ and we can obtain
(see \cite {ALM}, Lemma 5.5):

\begin{lemma} \label {geometric puzzle}

Let $0<\phi<\psi<\gamma<\pi/2$ be fixed.  For arbitrarily big $k>0$,
if $|I_\n|/|I_{\n-1}|$ is small enough, there exists a sequence
$V^j$ of open Jordan disks such that
$D_\phi(A^j) \subset V^j \subset D_\psi(A^j)$ and
$V^0=D_{(\phi+\psi)/2}(A^0)$ with the following
properties:

(1)\, If $j \neq 0$ and
$f(A^j) \subset A^k$ then $f:V^j \to V^k$ is a diffeomorphism;

(2)\, If $f(A^0) \cap A^j \neq \emptyset$,
then $\mod f(V^0) \setminus \overline {D_\gamma(A^j)}>k$.

\end{lemma}

\subsection{A special Banach space of perturbations}

Let $A^1=[l,r]$ with $l<r$, and let $N=[-l,l]$.
Domains $V^j$ which do not intersect $A^1$ or $N$
will play no role in the construction to follow.
Let $V$ be the union of all $V^j$ such that $A^j \subset N \cup A^1$.

One of the main problems of \cite {ALM} is to obtain a direction $v$ (or
infinitesimal perturbation) which
is transverse to the topological class of $f$.  The idea is to consider a
perturbation which does not affect much $f$ in $N$, but causes a bump near
the critical value, localized in $A^1$.  There are several difficulties
related to this scheme, the first of which is that such a bump can only be
reasonably controlled up to its first derivative.  Another difficulty is
that we want an analytic perturbation, so it cannot vanish in $N$ and be a
bump at $A^1$.  The solution involves the consideration of
certain Banach spaces of smooth ($C^1$) functions in $N \cup A^1$ which are
analytic in $\inter N \cup \inter A^1$,
which allows one to construct perturbations that, while badly behaved in the
real line (can be only controlled up to the first derivative), are well
behaved with respect to the complex puzzle structure.

While the proof in \cite {ALM} involves two steps, construction of a
transverse smooth vector field and approximation of this vector field
by polynomials, which need two different Banach spaces, we will realize
the same construction with just one Banach space.  This is important
to estimate the asymmetric roles of perturbations
concentrated in $N$ and $A^1$.  The proof of our main perturbation
estimate (Lemma \ref {puzzleest}) is a mixture of two estimates,
Lemma 7.4 (for perturbations localized in $A^1$) and Lemma 7.9
(for perturbations supported on $N \cup A^1$) of \cite {ALM}.

Let $Z=D_\gamma(A^1) \cup D_\gamma(N)$, and let $\Upsilon$ be the space of
all vector fields $v$ holomorphic on $Z$ and whose derivative admits
a continuous extension to $\overline Z$, which vanish up to the first
derivative in $\partial A^1$ and its forward iterates (this is a finite
set) and such that $v|D_\gamma(N)$ is a symmetric (odd) vector field.  We
use the norm $\|v\|=\sup_{z \in Z} |D v|$.

Let $\Upsilon=\Upsilon_1 \oplus \Upsilon_2$, where $v \in \Upsilon_1$ if
$v|D_\gamma(N)=0$ and $v \in \Upsilon_2$ if $v|D_\gamma(A^1)=0$.

Let $f_v=f \circ (\id +v)$.
The reader should think of vector fields $v \in \Upsilon$ as perturbations
of $f$ acting by $v \to f_v$.
One of the main advantages of the definition of
$\Upsilon$ is that, for small $v \in \Upsilon$, ``the puzzle persists'',
that is, there exists a continuation $V^v$ of the set $V$ inside $Z$, whose
connected components behave, under iteration by $f_v$, in the same way that
the connected components of $V$ behaved under iteration by $f$.

To make this more precise, let us say that $v \in \Upsilon$ is admissible if
there exists a holomorphic motion $h^v$ over $\D$, which is real-symmetric
if $v$ is real-symmetric, and is defined by the family of
transition maps $h^v[0,\lambda] \equiv h^v_\lambda:\C \to \C$,
$\lambda \in \D$ such that:

\begin{enumerate}

\item $h^v_\lambda|\C \setminus Z=\id$,
$h^v_\lambda|\partial f(V^0)=\id$;

\item $f_{\lambda v} \circ h^v_\lambda|V \setminus V^0=
h_\lambda \circ f$, $f_{\lambda v} \circ h^v_\lambda|\partial V^0=f$.

\end{enumerate}

The holomorphic motion $h^v$ will be said to be {\it compatible}
with $v$.

The following is a restatement of Lemma 7.9 of \cite {ALM}.

\begin{lemma} \label {7.9}

There exists $\epsilon>0$ such that if $v$ belongs to $\{v \in
\Upsilon|\|v\|<\epsilon\}$ then $v$ is admissible.

\end{lemma}

We also need the following simple estimate (see the proof of Lemma 7.4 of
\cite {ALM}):

\begin{lemma} \label {7.4}

Let $0<\theta<\gamma<\pi/2$.  There exists $\epsilon'>0$ such that if
$A$ is an interval and
$v$ is holomorphic on $D_\gamma(A)$ whose derivative extends
continuously to $\overline {D_\gamma(A)}$ satisfying $|D v|<\epsilon'$
then $\id+v:D_\gamma(A) \to \C$ is a
diffeomorphism and $D_\theta(A) \subset (\id+v) (D_\gamma(A))$.

\end{lemma}

Now we can prove:

\begin{lemma} \label {puzzleest}

There exists constants $\epsilon_1>0$, $\epsilon_2>0$, {\it where
$\epsilon_1$ depends only on $\psi$ and $\gamma$}
such that if $v_1 \in \Upsilon_1$, $\|v_1\|<\epsilon_1$ and
$v_2 \in \Upsilon_2$, $\|v_2\|<\epsilon_2$ then $v=v_1+v_2$ is admissible.

\end{lemma}

\begin{pf}

Let $n_1$ be such that $f^{n_1}(V^1)=V^0$
and let $\theta=(\psi+\gamma)/2$.

Let $v \in \Upsilon$ with $\|v\|<\epsilon$.  By Lemma \ref {7.9},
there exists a holomorphic motion $h^v$ compatible with $v$.

We claim that if $0<\epsilon_2<\epsilon$ is small enough
and $\|v\|<\epsilon_2$ then
for $\lambda \in \D$, $h^v_\lambda(V^1) \subset D_\theta(A^1)$.
Indeed, if this is not the case,
there would be a sequence $z_k \in \partial D_\theta(A^1)$,
$v_k \in \Upsilon$, $v_k \to 0$, such that
$f_{v_k}^{n_1+1}(z_k) \in f(V^0)$.
It clearly follows that $z_k \to \partial A^1=\{l,r\}$, let us say that $z_k
\to l$.  It is clear that
$$
f_{v_k}^{n_1+1}(z_k)=
f_{v_k}^{n_1+1}(l)+Df_{v_k}^{n_1+1}(l) z_k+o(z_k)=
f^{n_1+1}(l)+Df^{n_1+1}(l) z_k+o(z_k).
$$
In particular, the sequence $f_{v_k}^{n_1+1}(z_k)$ converges to
$f^{n_1+1}(l)$ along a direction which makes angle $\theta$ with the
real line (since $Df^{n_1+1}(l) \in \R \setminus \{0\}$), so
$f_{v_k}^{n_1+1}(z_k) \notin f(V^0)$ for $k$ big, which is a contradiction.

Let $\epsilon_1$ be as in Lemma \ref {7.4}.
If $v=v_1+v_2$, with $v_i \in \Upsilon_i$
and $\|v_i\|<\epsilon_i$, let $h^v_\lambda:\C \setminus (D_\gamma(A^1)
\setminus V^1)$ be given by $h^v_\lambda|(\C \setminus
D_\gamma(A^1))=h^{v_2}_\lambda$ and
$h^v_\lambda|V^1=((\id+\lambda_{v_1})|D_\gamma(A^1))^{-1}
\circ h^{v_2}_\lambda$.  Any extension of $h^v_\lambda$ to $\C$ is clearly
compatible with $v$.
\end{pf}

We will also need the following easy lemma:

\begin{lemma} \label {small central}

If $|I_\n|/|I_{\n-1}|$ is sufficiently small, then for $w=w_1+w_2$ with
$w_i \in \Upsilon_i$, $\|w_i\|<\epsilon_i$,
and for $\lambda \in \D$, then
$(f_{\lambda w}|h^w_\lambda(V^0))^{-1} (D_\gamma(V^1))
\subset \D_{\rho |A^0|}(0)$,
where $\rho \to 0$ as $|I_\n|/|I_{\n-1}| \to 0$.

\end{lemma}

\begin{pf}

Let $U=h_\lambda(V^0)$ and
$U^0=(f_{\lambda w}|W)^{-1} (D_\gamma(A^1))$.  Notice that $f_{\lambda
w}(0)=f(0) \in D_\gamma(V^1)$.  Thus, $f_{\lambda w}|(U \setminus \overline
{U^0})$ is a double covering of $f(U_0) \setminus \overline
{D_\gamma(A^1)}$.  By Lemma \ref {geometric puzzle},
if $|I_\n|/|I_{\n-1}|$ is small then
$\mod(f(U_0) \setminus \overline {D_\gamma(A^1)})$ is large, and so
$\mod (U \setminus \overline {U^0})$ is also big.  Since the derivative of
$\id+\lambda w$ is smaller than $\max \{1+\epsilon_1,1+\epsilon_2\}$,
we see that the diameter of $U$ is at most $2 |A^0|$, so
the diameter of $U^0$ can be bounded by $\rho |A^0|/2$ with small $\rho$
as required.
\end{pf}

\subsection{Analytic vector fields}

We will be specially concerned with special types of $w$ which generate
analytic families of unimodal maps.  The following lemma is obvious:

\begin{lemma} \label {6.6}

If $w \in \Upsilon$ is real-symmetric and
has an analytic extension $w:I \to I$ of $C^1$ of norm less than one,
such that $w(-1)=w(1)=0$, then $f_{\lambda w}$,
$\lambda \in (-1,1)$ is an analytic family of unimodal maps, and $I_\n$ is a
nice interval with preperiodic boundary for each $f_{\lambda w}$.

\end{lemma}

The following is a consequence of the Mergelyan Polynomial Approximation
theorem:

\begin{lemma}

Let $w \in \Upsilon$.  Then there exists a sequence
$w_m \in \Upsilon$ such
that the $C^1$ norm of $w_m|I$ is less than $\|w\|$,
$w_m(-1)=w_m(1)=0$ and $w_m \to w$ in $\Upsilon$.  If $w$ is real-symmetric
then we can also choose $w_m$ real-symmetric.

\end{lemma}


\begin{lemma}

Let $w \in \Upsilon$ be as in Lemma \ref {6.6}.  If $w$ is admissible,
then the domain of the first return map to $I_\n$ under iteration by
$f_{\lambda w}$ is
$((\id+\lambda w)|h^w_\lambda(V^0))^{-1}(V) \cap \R$.

\end{lemma}

\begin{pf}

By construction,
all components of $((\id+\lambda w)|h^w_\lambda(V^0))^{-1}(V) \cap \R$
are components of the first return map
to $I_\n$, so we just have to check that all components are of this form.
Notice that each $x \in V \cap (f(-l),l)$ has two preimages by $f$
in $V \cap ((-l,l) \setminus I_\n)$.  It follows that each
$x \in h^w_\lambda(V) \cap (f(-l),l)$ has two preimages by $f^w_\lambda$
in $h^w_\lambda (V) \cap ((-l,l) \setminus I_\n)$.
Let now $T$ be a component
of the first return map to $I_\n$ under iteration by $f_{\lambda w}$.
If $T$ is the
central component, then $T$ must be the preimage of $A_1$. 
Otherwise, all iterates of $T$ up to the return are contained in 
$(f(-l),l)$.  Since $\inter I_\n \subset h^w_\lambda(V)$, we conclude that
all iterates of $T$ up to the return belong to $h^w_\lambda(V)$.
\end{pf}

\comm{
Let us consider a component $T_0$ of the first return map to $I_\n$
under $f_{\lambda w}$, and let $T_i=f_{\lambda w}^i(T_0)$, $i=1,...,l$,
where $l$ is the time of first return so that $T_l \subset I_\n$.
Notice that $T_i$, $1 \leq i <l$ are contained in
$A_1 \cup N$ (which is forward invariant).  Notice that all $z \in V \cap
\R$

\item there exists homeomorphisms $H[\lambda]$,
depending continuously on $\lambda$ such that $f_{\lambda w} \circ
H[\lambda]|I \setminus I_\n=H[\lambda] \circ f$,

\item $H[\lambda](I_\n)=I_\n$,

\end{enumerate}

\end{lemma}
}

\subsection{A special perturbation}

Let us consider an affine map $Q:A^1 \to I$, and let
$$
\tilde v_n(z)=(1-z^2)(1-e^{-2n})+
\frac {2} {n}(e^{-n(1+z)}+e^{-n(1-z)}-e^{-2n}-1),
$$
and let $v_n \in \Upsilon_1$ be such that
$v_n|D_\gamma(A^1)=Q^* \tilde v_n \epsilon_1/8$.
Notice that $\|v_n\|<\epsilon_1$.

\subsubsection{Infinitesimal transversality}

The importance of the sequence $v_m$ in \cite {ALM} is that it is eventually
transverse to the topologically class of $f$.

\comm{
A quasiconformal vector field $\alpha$ of $\overline \C$ is a continuous
vector field with locally integrable distributional derivatives
$\op \alpha$ and $\partial \alpha$
in $L^1$ and $\op \alpha \in L^\infty$.
}

Let us say that $w$ is {\it formally transverse at $f$} if there is no
quasiconformal vector field $\alpha$ of $\C$, such that
for $z \in \orb_f(0)$, $w(z)=f^* \alpha(z)-\alpha(z)$.
(This definition is motivated by Theorem \ref {8.4}, see also Lemma \ref
{transv}.)

The following summarizes Lemmas 7.6, 7.7 and 7.8 of \cite {ALM}.

\begin{lemma} \label {v_m}

Let $v_m$ be defined as above.  If $|I_\n|/|I_{\n-1}|$ is sufficiently
small, then for $m$ sufficiently big, $v_m$ is formally transverse at $f$.

\end{lemma}

\comm{
\begin{lemma} \label {transv}

If $w$ is an analytic vector field on $I$ and there exists no quasiconformal
vector field $\alpha$ satisfying $w=f^* \alpha-\alpha$ on the critical
orbit, then the analytic family of
unimodal maps $f_{\lambda w}$, $\lambda \in (-\epsilon,\epsilon)$ is
transverse to the topological class of $f$.

\end{lemma}

\begin{pf}

Indeed, we have
$$
\frac {d} {d\lambda} f_{\lambda v}|_{\lambda=0}=Df w.
$$
By hypothesis, there exists no qc vector field $\alpha$ satisfying $Df
w=\alpha \circ f-\alpha Df$ on the critical orbit.  By Theorem 8.10 of \cite
{ALM}, $Df w$ is transverse to hybrid class of $f$.
\end{pf}
}

The following is due to (a version of) the so called Key estimate of
\cite {ALM} (more precisely we use Corollary 7.14 of \cite {ALM}):

\begin{lemma} \label {w closed}

The set of vector fields $w \in \Upsilon$ which are not formally transverse
at $f$ is a closed subspace of $\Upsilon$.

\end{lemma}

\begin{rem} \label {w form}

In particular, if $m$ is sufficiently big and $w$ is close to $v_m$ then $w$
is formally transverse at $f$.
\end{rem}

\subsubsection{Macroscopic transversality}

The following result can be interpreted as the
macroscopic counterpart to the infinitesimal transversality of $v_m$.

Let $r>0$ be minimal with $f^{r+1}(0) \in V^1$.

\begin{lemma} \label {image}

There exists a constant $\tau_0>0$ depending only on
$\epsilon_1$ and $\phi$, such that if $|I_\n|/|I_{\n-1}|$ is sufficiently
small the following holds.  Let $v_m$ be defined as above and let
$r>0$ be minimal with $f^{r+1}(0) \in V^1$.
Then for $m$ sufficiently big,
there exists a domain $\hat \Theta \subset \D$ such that the map
$\theta:\hat \Theta \to \C$ given by $\theta(\lambda)=f^r_{\lambda v_m}(0)$
is a diffeomorphism onto $\D_{\tau_0 |I_\n|}$.

\end{lemma}

\begin{pf}

Since $\|v_m\|<\epsilon_1$,
there exists a holomorphic motion $h^{v_m}$
which is compatible with $v_m$.

Let $\Psi:\D \to \C$, $\Psi(\lambda)=(\id+\lambda v_m)(f(0))$.  It is
clearly a diffeomorphism over a round disk $D_m$ centered on $0$.
Let $d_m$ be the hyperbolic distance between $f(0)$ and $\partial D_m$ in
$D_{\pi/2}(A^1)$.  It is easy to estimate directly
$d_m$ from below in terms of $\epsilon_1$ and $m$.  In particular, for
$m$ big, $d_m>\tilde \tau>0$ where $\tilde \tau$ depends only on
$\epsilon_1$, not on the position of $f(0)$ in $A^1$.\footnote {To see
this, notice that $D\Psi(0)=v_m(f(0))$, and the norm of $v_m(f(0))$ in the
hyperbolic metric of $D_{\pi/2}(A^1)$ at $f(0)$ is at least $\epsilon_1/10$
for $m$ big.  Let $P:D_{\pi/2}(A^1) \to \D$ be a Moebius transformation
taking $f(0)$ to $0$.  The the norm of $D(P \circ \Psi)(0)$ in the
hyperbolic metric of $\D$ at $0$ is at least $\epsilon_1/10$, so the
Euclidean norm of $D(P \circ \Psi)(0)$ is at least $\epsilon_1/10$.  By the
Koebe 1/4 Theorem, $P(D_m)$ contains a round disk of radius
$\epsilon_1/40$ around $0$, thus the hyperbolic distance from
$\partial P(D_m)$ to $0$ in $\D$ is at least $\epsilon_1/40$.}

Now let $Q$ be the connected component of $f(0)$ on
$f^{-(r-1)}(V^0)$, so that $f^{r-1}:Q \to V^0$ is a diffeomorphism.
The hyperbolic distance between $\partial D \cap Q$ and $f(0)$ in $Q$ is
bounded from below by $\tilde \tau$ by the Schwarz Lemma (if $\partial D
\cap Q=\emptyset$, we let this distance be $\infty$).  It follows that
$f^{r-1}(Q \cap D)$ contains a $\tilde \tau$ hyperbolic neighborhood of
$f^r(0)$ on $V^0$.
Now, if $|I_\n|/|I_{\n-1}|$ is very small, then $|I_{\n+1}|/|I_\n|$
is also very small, so $f^r(0)$ (which is contained in $I_{\n+1}$)
is $\tilde \tau/2$ close to $0$ in the hyperbolic metric of
$V^0 \supset D_\phi(A^0)$.

As a consequence, $f^{r-1}(Q \cap D)$
contains a $\tilde \tau/2$ hyperbolic neighborhood of $0$ in $V^0$, and
since $V^0 \supset D_\phi(A^0)$, it must contain $\D_{\tau |A^0|}$, where
$\tau$ depends on $\epsilon_1$ and $\phi$.
\end{pf}

\subsubsection{Construction of a full $R$-family}

Let $\tau_0$ be the constant of Lemma \ref {image} and
let $|I_\n|/|I_{\n-1}|$ be such that Lemma \ref {small central}
holds with $\rho<\tau_0/4$.

Let $m$ be big and let us fix $v=v_m$ such that
Lemmas \ref {image} and \ref {v_m} are
valid, and let $\hat \Theta$ be the domain of Lemma \ref {image}.

Let $w=w_1+w_2$ with $w_i \in \Upsilon_i$, $\|w_i\|<\epsilon_i$.

Let $U[0]=V^0$ and let the family $\{U^j[0]\}$ denote the connected
components of $(f|V^0)^{-1}(\cup V^j)$, letting $0 \in U^0[0]$.

Let us consider a holomorphic motion $\tilde H$ over $\D$ given by the
transition maps $\tilde H[0,\lambda]=\tilde H_\lambda:\C \to \C$ such that:
$$
\tilde H_\lambda|\C \setminus U[0]=h_{\lambda w}
$$
$$
f_{\lambda w} \circ \tilde H_\lambda|
U[0] \setminus U^0[0]=h_{\lambda w} \circ f.
$$

Let $U[\lambda]=\tilde H_\lambda
(U[0])$, $U^j[\lambda]=\tilde H_\lambda(U^j[0])$.

Let $R[\lambda]$ be the first return map from $U^j[\lambda]$ to $U_0$.
It is clear that $(R[\lambda],\tilde H_\lambda)$ has a
structure of a (non-full) $R$-family over $\D$.
Let us consider the landing family $(L[\lambda],H_\lambda)$ associated to
$(R[\lambda],\tilde H_\lambda)$.

Let $W^\d[0]$ be the domain of $L[0]$ containing $R[0](0)$.
Notice that $L[\lambda]|W^\d[\lambda]$
extends to a diffeomorphism $R^\d[\lambda]$ onto $U[\lambda]$.  For
$\tau<\tau_0$, let $\Delta_\tau[\lambda]$ be the
preimage of $\D_{\tau |A^0|}(0)$ by this diffeomorphism.

If $w=v$ then $R^\d[\lambda]=R^\d[0]$ for all
$\lambda$, since $v$ is supported on $D_\gamma(A^1)$.

In particular, $R^\d[\lambda]=R^\d[0]$ and
$\Delta_\tau[\lambda]=\Delta_\tau[0]$ for all
$\lambda$.  So $\lambda \mapsto R[\lambda](0)$ is a map
which restricts (in some domain $0 \in \O^v$)
to a diffeomorphism onto $\Delta_\tau[0]$.
It follows that taking $\tau=\tau_0/2$, for any $w$ close to $v$
there exists a domain $0 \in \O^w$ where
$\lambda \mapsto R[\lambda](0)$ is a diffeomorphism onto
$\Delta_\tau[0]$ (of course, $\O^w$ depends on
$w$).

But for $w \in \Upsilon$ close to $v$ and
for all $\lambda \in \D$, $U^0[\lambda]$ is contained in
$\D_{\rho |A^0|}$, so $W^\d[\lambda]$ is contained in
$\Delta_{\tau_0/2}[0]$ with space.
By the argument principle,
letting $\Theta$ be the connected component of $0$ of the set
of $\lambda \in \hat \Theta$ with $R[\lambda](0) \in W^\d[\lambda]$, the map
$S:\overline \Theta \to \overline {W^\d[0]}$ such that
$S(\lambda)=H_\lambda^{-1}(R[\lambda](0))$
is a homeomorphism.  We also have that the diameter of $\Theta$
is very small if
$\rho$ is small (in particular if $|I_\n|/|I_{\n-1}|$ is small).

Let $U_1[0]=U^0[0]$ and
let $\{U^j_1[0]\}$ be the connected components of the preimage by
$R[0]|U^0[0]$ of $\cup W^\d[0]$, and let $0 \in U^0_1$.

Let $h$ be a holomorphic motion over $\Theta$ given by transition maps
$h[0,\lambda]=h_\lambda:\C \to \C$ such that
$$
h_\lambda|\C \setminus U_1=H_\lambda,
$$
$$
R[\lambda] \circ h_\lambda|U_1 \setminus
U^0_1=h_\lambda \circ R[0].
$$

Let $U_1[\lambda]=h_\lambda(U_1[0])$ and
$U^j_1[\lambda]=h_\lambda(U^j_1[0])$.

Our construction shows clearly that the first return map
$R_1[\lambda]$ from $\cup U^j_1[\lambda]$ to
$U_1[\lambda]$ is an $R$-map for $\lambda \in \Theta$,
so $(R_1[\lambda],h_\lambda)$ is an $R$-family, and our choice of
$\Theta$ implies that $R_1[\lambda]$ is a {\it full} $R$-family.

Let us summarize the properties we obtained in this construction:

\begin{lemma} \label {special}

If $|I_\n|/|I_{\n-1}|$ is small enough, there exists a real-symmetric
vector field $v \in \Upsilon$ and a neighborhood
$v \in \VV \subset \Upsilon$ such that for any
$w \in \VV$ real-symmetric, there exists a
domain $0 \in \Theta \subset \D$, a family of $R$-maps
$R_1[\lambda]:U^j_1[\lambda] \to U_1[\lambda]$,
$\lambda \in \Theta$, and a real-symmetric holomorphic motion $h$ over
$\Theta$
such that:

(1)\, For $\lambda \in \Theta \cap \R$, $U_1[\lambda] \cap \R=I_{\n+1}$;

(2)\, $R_1[\lambda]$ is the first return map from $\cup U^j_1[\lambda]$ to
$U_1[\lambda]$ under iteration by $f_{\lambda w}$;

(3)\, $(R_1[\lambda],h)$ form a full real-symmetric $R$-family.

And moreover, if $w$ is as in Lemma \ref {6.6} and $\lambda \in \Theta \cap
\R$ then:

(4)\, $I_{\n+1}[\lambda] \equiv U_1[\lambda] \cap \R$ is the component of
$0$ of the first return map to $I_\n$ under iteration by $f_{\lambda w}$;

(5)\, $I^j_{\n+1}[\lambda] \equiv U^j_1[\lambda] \cap \R$ are the
domains of the first return map to $I_{\n+1}[\lambda]$ under iteration of
the real analytic extension $f_{\lambda w}:I \to I$.

\end{lemma}

The construction of the $R$-family gives us also a good control of its
geometry.

\begin{lemma} \label {special1}

In the setting of Lemma \ref {special}, $\Dil(h|\C \setminus
\overline {U^0_1})<1+\epsilon$, and $\mod (U_1[0] \setminus
\overline {U^0_1[0]})>C$, where $\epsilon \to 0$ and $C \to \infty$ when
$|I_\n|/|I_{\n-1}| \to 0$.

\end{lemma}

\begin{pf}

Indeed, $\Dil(h|\C \setminus
\overline {U^0_1})<1+\epsilon$ is bounded by the hyperbolic diameter
of $\Theta$ in $\D$, which is small if $|I_\n|/|I_{\n-1}| \to 0$ is big.
On the other hand, $\mod (U_1[0] \setminus
\overline {U^0_1[0]}) \geq \mod (U[0] \setminus
\overline {U^0[0]})/2 \geq \mod (f(V^0) \setminus \overline {V^1})/4>k/4$,
which is big if $I_\n \setminus I_{\n-1}$ is small by Lemma \ref {geometric
puzzle}.
\end{pf}

\subsection{Remarks on the infinitesimal transversality of the special
perturbation}

We would like to point out that the ``macroscopic transversality'' of
$v_m$ is very much related to its infinitesimal transversality.  The
(formalizable) argument relating both properties is as follows (notice that
this argument is different from the one given in \cite {ALM}, which
emphasizes estimates at the infinitesimal level):

(1) $v_m$ can be $C^1$ extended to $I$ as an odd vector field
which vanishes on $[r,1]$, $[-1,-r]$ and $[-l,l]$.  This vector field is not
$C^2$ but its $C^1$ norm is small ($\epsilon_1$).

(2) (Macroscopic transversality implies a $C^1$ connecting lemma)
Notice that the interval
$(f_{-v_m}^r(0),f_{v_m}^r(0))$ contains the interval $I_{\n+1}$
(with lots of space).  We conclude that
the family $f_{\lambda v_m}$, $\lambda \in (-1,1)$
must go through a parameter $\lambda$ where $f_{\lambda v_m}^r(0)=0$, and so
changes the combinatorics of $f$.

(3) Using the Key Estimate of \cite {ALM}, we see that, if $v_m$ is not
formally transverse at $f$, then it is actually tangent to the topological
class of $f$ in the following sense.  There exists a (real-symmetric)
holomorphic motion $h$ over $\D$ whose transition maps
$h[0,\lambda] \equiv h_\lambda:\C \to \C$ are such that
$f_\lambda=h_\lambda \circ f \circ h_\lambda^{-1}$
is a family of so called ``puzzle maps'' (which behave as unimodal maps)
such that
$$
\left .\frac {d} {d\lambda} f_\lambda \right |_{\lambda=0}=
\left .\frac {d} {d\lambda} f_{\lambda v_m} \right |_{\lambda=0}=
Df \cdot v_m
$$
(the maps $h_\lambda$ are characterized by $\op h_\lambda/\partial
h_\lambda=\lambda \partial \alpha$ for a specially chosen
quasiconformal vector field $\alpha$ satisfying
$v_m=f^* \alpha-\alpha$ on the critical orbit).  This family can be
considered the Beltrami path through $f$ in the direction of
$Df \cdot v_m$.

(4) The family $f_\lambda$ is tangent to $f_{\lambda v_m}$ at $\lambda=0$
and both families have big extensions (to $\D$).  In particular, they must
be close together in a slightly smaller disk, where we can detect
the change of
combinatorics: there is a parameter $\lambda \in \D$ such that
$f_\lambda^r(0)=0$\footnote
{More precisely, we use that the holomorphic map
$\lambda \mapsto f^r_\lambda(0)$ has the same
derivative at $0$ as the almost linear map
$\lambda \mapsto f^r_{\lambda v_m}(0)$, and a simple estimate shows
that there exists a parameter $\lambda \in \D$ such that $f^r_{\lambda
v_m}(0)=0$.}.

(5) In particular, the family $f_\lambda$ must change combinatorics, but
this is a contradiction, since it consists of maps topologically conjugate
to $f$.  So we conclude that $v_m$ is formally transverse at $f$.
Notice that our argument is that a ``reasonably efficient''\footnote{In
the sense of admitting a controlled extension to a big domain, as the
Beltrami path we constructed.} tangent path to $v_m$ closes
macroscopically the critical orbit.

(6) (Infinitesimal analytic connecting lemma)
Although $v_m$ is only $C^1$ in the interval,
we can approximate it in the topology of $\Upsilon$
by polynomials $w$ which
will be still formally transverse to $f$.  Those vector fields $w$ are
transversal to the topological class of $f$:
they close ``infinitesimally'' the critical orbit.

\section{The Phase-Parameter relation} \label {phase-parameter}

\subsection{Phase-Parameter relation for the special family}

Let $f \in \FF$ and let
$R_i:\cup I^j_i \to I_i$ be the first return map.  For $\d \in \Omega$,
$\d=(j_1,...,j_m)$, let $I^\d_i=\{x \in I_i|\,
R^{k-1}_i(x) \in I^{j_{i+1}}_i, 1 \leq k \leq m\}$, and let
$R^\d_i=R_i^m|I^\d_i$.  Let
$C^\d_i=(R^\d_i)^{-1}(I^0_i)$.  The map $L_i:\cup C^\d_i \to I^0_i$ is the
first landing map from $I_i$ to $I_{i+1}$.

\begin{definition}

Let us say that a family $f_\lambda$ of unimodal maps satisfies the
Topological Phase-Parameter relation at a parameter $\lambda_0$ if
$f=f_{\lambda_0} \in \FF$, and there exists $i_0>0$ and
a sequence of nested intervals $J_i$, $i \geq i_0$ such that: 

\begin{enumerate}

\item $J_i$ is the maximal interval containing $\lambda_0$ such that for
all $\lambda \in J_i$ there exists a homeomorphism $H_i[\lambda]:I \to I$
such that
$f_\lambda \circ H_i[\lambda]|(I \setminus I_{i+1})=H_i[\lambda] \circ f$.

\item There exists a homeomorphism $\Xi_i:I_i \to J_i$ such that
$\Xi_i(C^\d_i)$ (respectively, $\Xi_i(I^\d_i)$)
is the set of $\lambda$ such that
the first return of $0$ to $H_i[\lambda](I_i)$
under iteration by $f_\lambda$ belongs to $H_i[\lambda](C^\d_i)$
(respectively, $H_i[\lambda](I^\d_i)$).

\end{enumerate}

\end{definition}

\begin{definition}

Let $f_\lambda$ be a family of unimodal maps.  We say that $f_\lambda$
has Decay of Parameter Geometry at $\lambda_0$ if $f=f_{\lambda_0} \in \FF$,
it satisfies the Topological Phase-Parameter relation at $\lambda_0$ and
$|J_{n_k+1}|/|J_{n_k}|<C \lambda^k$ for some constants $C>0$, $\lambda<1$,
where $n_k-1$ counts the non-central levels of the principal nest of $f$.

\end{definition}

\begin{thm} \label {7.1}

Let $f \in \FF$ be analytic.  There exists a polynomial vector field $w$
such that the family $f_{\lambda w}=f \circ (\id+\lambda w)$,
$\lambda \in (-\epsilon,\epsilon)$ is an analytic family of unimodal maps
which satisfies the Topological Phase-Parameter relation and has
Decay of Parameter Geometry at $0$.

\end{thm}

\begin{pf}

Let $w$ and $\n$ be as in Lemma \ref {special}.
Denote by $(\RR_1,h_1)$ the $R$-family of that
lemma.  Since $f \in \FF$, the critical point is recurrent and we can
clearly construct a $R$-chain $(\RR_i,h_i)$ over $\lambda=0$.  It is clear
that the real trace of $R_i[0]:\cup U^j_i[0] \to U_i[0]$ is the first return
map to $I_{\n+i}$.  Let $J_{\n+i}=\Lambda_i \cap \R$, let
$\Xi_{\n+i}=\chi_i[0]|I_{\n+i}$.  It is clear that $|J_{n_k+1}|/|J_{n_k}|$
decays exponentially by Lemma \ref {special1} and Theorem \ref {3.8},
where $n_k-1$ counts the non-central levels
of the principal nest of $f$.  In particular, $|J_n| \to 0$.

In order to conclude the result, we just have to show the existence of the
continuous family of homeomorphisms
$H_i[\lambda]$, for $i$ sufficiently big.  Notice that if $\lambda \in
J_{\n+i}$, if $p \in I_{\n+i}[\lambda]$ is a periodic orbit for
$f_\lambda$ which never enters $I^0_{\n+i}[\lambda]$
then $p$ is hyperbolic by the Schwarz Lemma.  So, if $\lambda \in J_{\n+i}$,
the only non-hyperbolic periodic orbits for $f_\lambda$ must be entirely
contained in $I \setminus I_{\n+1}$.  But since $f|I \setminus I_{\n+1}$ is
hyperbolic, there exists $\epsilon>0$ such that if $\lambda \in
(-\epsilon,\epsilon)$, all periodic orbits in $I \setminus
I_{\n+1}[\lambda]$ of $f_\lambda$ are hyperbolic (by Lemma \ref {5.2}).
In particular, if $i$ is
sufficiently big, $J_i \subset (-\epsilon,\epsilon)$, and all periodic
orbits of $f_\lambda$ in $I \setminus I_{i+1}[\lambda]$ are hyperbolic.
The result follows by Lemma \ref {5.3}.
\end{pf}

Let $K_i$ be the closure of the union of all
$\partial C^\d_i$ and $\partial I^\d_i$.
Notice that $H_i$ and $\Xi_i$ are only uniquely defined in $K_i$.
Condition (2) of the Topological Phase-Parameter relation can be
equivalently formulated as the existence of a homeomorphism $\Xi_i:I_i \to
J_i$ such that the first return of the critical point (under iteration by
$f_\lambda$) to $H_i[\lambda](I_i)$ belongs to
$H_i[\lambda](K_i)$ if and only if $\lambda \in \Xi_i(K_i)$.

Let us now estimate the metric properties of $H_i|K_i$ and $\Xi_i|K_i$.  In
order to do so, we will need to consider convenient restrictions of those
maps.

Let $\tilde I_{i+2}=(R_i|I^0_i)^{-1}(I^\d_i)$, where $\d$ is such that
$(R_i|I^0_i)^{-1}(C^\d_i)=I_{i+2}$.

Let $\tau_i$ be such that $R_i(0) \in I^{\tau_i}_i$.

Let $\tilde K_i=\overline {(\cup_j \partial I^j_i \cup \partial I_i)
\setminus \inter \tilde I_{i+1}}$.

Let $J^j_i=\Xi_i(I^j_i)$.

\begin{definition}

Let $f_\lambda$ be a family of unimodal maps.  We say that $f_\lambda$
satisfies the Phase-Parameter relation at $\lambda_0$ if $f=f_{\lambda_0}$
is simple, $f_\lambda$ satisfies the
Topological Phase-Parameter relation at $\lambda_0$ and for every $\g>1$,
there exists $i_0>0$ such that for $i>i_0$ we have:

\begin{description}

\item[PhPa1] $\Xi_i|(K_i \cap I^{\tau_i}_i)$ is $\g$-qs,

\item[PhPa2] $\Xi_i|\tilde K_i$ is $\g$-qs,

\item[PhPh1] $H_i[\lambda]|K_i$ is $\g$-qs if $\lambda \in
J^{\tau_i}_i$,

\item[PhPh2] the map $H_i[\lambda]|\tilde K_i$ is
$\g$-qs if $\lambda \in J_i$.

\end{description}

\end{definition}

\begin{thm} \label {7.2}

In the same setting of the previous theorem, if $f$ is simple, the family
$f_{\lambda w}$ satisfies the Phase-Parameter relation at $0$.

\end{thm}

\begin{pf}

Let $(\RR_i,h_i)$ be the $R$-chain of the proof of Theorem~\ref {7.1}.
By Theorems \ref {3.4} and \ref {3.8},
$\mod(U_i[0] \setminus \overline {U^0_i[0]}) \to \infty$ and
$\mod(\Lambda_i \setminus \overline {\Lambda_{i+1}}) \to \infty$
(notice that since $f$ is simple, all deep
enough levels $\RR_i$ are non-central).  This implies that, for any fixed
$\g>1$, there exists $i_0>0$ such that for $i>i_0$ the hypothesis of Lemmas
\ref {cphpa1} and \ref {cphpa2} are fulfilled and hence their conclusions
(CPhPa1, CPhPh1, CPhPa2, and CPhPh2) apply.  Those immediately imply the
four conditions (PhPa1, PhPh1, PhPa2, and PhPh2)
of the Phase-Parameter relation by restriction to the real line.
\end{pf}


\subsection{Phase-parameter relation for transverse families}

Let $f_{\lambda w}$ be the special family constructed before.

\begin{lemma} \label {transv}

The family $f_{\lambda w}$ is transverse to the topological class of $f$ at
$\lambda=0$.

\end{lemma}

\begin{pf}

Indeed, if $f_{\lambda w}$ is not transverse then by Theorem \ref {8.4},
there exists a qc vector field $\alpha:\C \to \C$ such that
$$
w Df=\left .\frac {d} {d\lambda} f_{\lambda w} \right |_{w=0}=
\alpha \circ f-\alpha Df
$$
on $\orb_f(0)$.
Dividing by $Df$ we get $w=f^* \alpha-\alpha$ on $\orb_f(0)$.  But this
contradicts Remark \ref {w form}.
\end{pf}

We will now show how to use the lamination of \cite {ALM} to transfer the
Phase-Parameter relation from the transversal family $f_{\lambda w}$ to
any transversal family $f_\lambda$.
The basic idea is contained in the following diagram:
$$
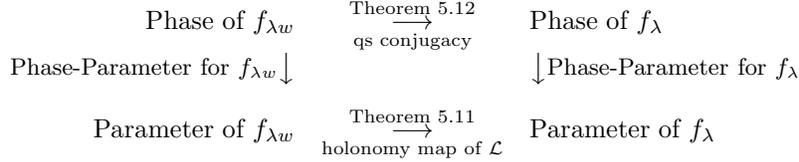
\begin{array}{rcl}
\text {Phase of $f_{\lambda w}$}& \underset {\text {qs conjugacy}}
{\overset {\text {Theorem \ref {8.2}}} {\longrightarrow}}
& \text {Phase of $f_\lambda$}\\
\text {\small Phase-Parameter for $f_{\lambda w}$}
\big\downarrow & & \big\downarrow \text {\small Phase-Parameter for
$f_{\lambda}$}\\[3mm]
\text {Parameter of $f_{\lambda w}$}&
\underset {\text {holonomy map of $\LL$}}
{\overset {\text {Theorem \ref {8.2'}}} {\longrightarrow}} &
\text {Parameter of $f_\lambda$}\\
\end{array}
$$
(notice that the estimates for all arrows are all ultimately based on
the $\lambda$-Lemma).

\begin{thm} \label {phpatra}

Let $f \in \FF$, and let $f_\lambda$ be a one-parameter
analytic family of unimodal maps
through $f$ such that $f_{\lambda_0}=f$ and $f_\lambda$ is transverse to the
topological class of $f$ at $\lambda=\lambda_0$.  Then the Topological
Phase-Parameter relation and Decay of Parameter Geometry holds for the
family $f_\lambda$ at $\lambda_0$.  Moreover, if $f$ is simple, then the
Phase-Parameter relation also holds.

\end{thm}

\begin{pf}

Using Theorems \ref {7.1}, \ref {7.2} and Lemma \ref {transv}
consider the family $f_{\lambda w}$ through $f$, which is
transverse to the hybrid class of $f$ and which satisfies the Topological
Phase-Parameter relation and Decay of Parameter Geometry (and the
Phase-Parameter relation if $f$ is simple).
Fix $a$ such that both $f_{\lambda w}$ and $f_\lambda$ are
analytic paths in $\U_a$.  Let $\LL$ be the lamination from
Theorem \ref {8.1}.  Since both $f_\lambda$ and
$f_{\lambda w}$ are transverse to the topological class of $f$ (at
$\lambda_0$ and $0$), we can consider the local holonomy map of
the lamination $\LL$, $\psi:(-\epsilon,\epsilon) \to
(\lambda_0-\epsilon',\lambda_0+\epsilon')$.

Let $\tilde \Xi_i:I_i \to \tilde J_i$ be the phase-parameter map for the
family $f_{\lambda w}$, and let $\tilde H_i[\lambda]$ be the phase-phase
map.  We obtain the phase-parameter map for $f_\lambda$
as a composition $\Xi_i=\psi \circ \tilde \Xi_i$.  Since $|\tilde J_i| \to
0$,
$$
\lim_{i \to \infty}
\sup_{\lambda \in \tilde J_i} \|f_{\lambda
w}-f_{\psi(\lambda)}\|_a=0.
$$
In particular, by Theorem \ref {8.2'},
$\psi|\tilde J_i$ is $\g_i$-qs with $\lim \g_i=1$.

Since for each $\lambda \in J_i=\psi(\tilde J_i)$, $f_\lambda$
is qs conjugate to $f_{\psi^{-1}(\lambda) w}$, we see that if
$\lambda \in J_i$ then there are no non-hyperbolic periodic orbits for
$f_\lambda$ in the complement of the continuation of $I_{i+1}$.
Using Lemma \ref {5.2} we conclude as in Theorem \ref {7.1}
the existence of a continuous family
$H_i[\lambda]$ of phase-phase maps for the family $f_\lambda$.
It follows that the Topological Phase-Parameter relation holds for
$f_\lambda$ at $\lambda_0$.

Since $\psi$ is quasisymmetric, it is H\"older and
the Decay of Parameter Geometry also follows from Theorem \ref {7.1}.
If $f$ is simple, estimates PhPa1 and PhPa2 follow from Theorem \ref {7.2}.

Let $h_\lambda:I \to I$ be a quasisymmetric conjugacy between
$f_{\lambda w}$ and $f_{\psi(\lambda)}$ which is $1+O(\|f_{\lambda
w}-f_{\psi(\lambda)}\|_a)$-qs.  This family might not be continuous, but
$H_i[\psi(\lambda)]|K_i=h_\lambda \circ \tilde H_i[\lambda]$,
which is enough for our purposes.  In particular,
if $f$ is simple, PhPh1 and PhPh2 follow from Theorem \ref {7.2}.
\end{pf}

\begin{rem}

Notice that even if we are only interested in the phase-parameter relation
for individual families, this proof needs the knowledge of the behavior of
topological conjugacy classes of unimodal maps
in infinite dimensional spaces.  For the case of the
quadratic family, this is not needed: the argument of \cite {parapuzzle} is
based on the combinatorial theory of the Mandelbrot set (Douady-Hubbard,
Yoccoz), which allows to show directly that the real quadratic family gives
rise to full unfolded complex return type families.  In particular,
our proof also gives a somewhat different approach to the
phase-parameter relation on the quadratic family itself.

\end{rem}

\section{Proof of Theorem A} \label {thmA}

Let $f_\lambda$ be a one-parameter non-trivial analytic family of
unimodal maps.  In view
of Theorem \ref {phpatra}, to conclude Theorem A it is enough to show that

\begin{enumerate}

\item Almost every non-regular parameter belongs to $\FF$, that is,
it is Kupka-Smale, has a recurrent critical point and
is not infinitely renormalizable,

\item Almost every parameter in $\FF$ is simple,

\item $f_\lambda$ is transverse to the topological class of almost every
parameter.

\end{enumerate}

We will take care of these issues separately below: item (1) will follow
from Lemmas \ref {ks}, \ref {nonre}, and \ref {infrenor}, item (2) from
Lemma \ref {simple} and item (3) from Lemma \ref {aetra}. 

\subsection{Transversality}

\begin{lemma} \label {ks}

Let $f_\lambda$ be a non-trivial analytic family of unimodal maps.  Then at
most countably many parameters are not Kupka-Smale or have a periodic or
preperiodic critical point.

\end{lemma}

\begin{pf}

Indeed, the set of parameters which are not Kupka-Smale correspond to
solutions of countably many analytic equations of the type
$f_\lambda^n(p)=p$, $Df_\lambda^{2n}(p)=1$, $n>0$.  Similarly, the set of
parameters with periodic or preperiodic critical point corresponds to
countably many equations of the type $f_\lambda^m(0)=f_\lambda^n(0)$,
$0 \leq m<n$.  So the set of parameters which
are not Kupka-Smale is either countable or contains intervals.  Since
regular parameters are dense, the first possibility holds.
\end{pf}

The following result is due to Douady, see Lemma 9.1 of \cite {ALM}:

\begin{lemma} \label {9.2}

Let $\LL$ be a codimension-one complex lamination on an open set $\VV$ of
some Banach space, and let $\gamma$ be an analytic path in $\VV$.  If
$\gamma$ is not contained in a leaf of $\LL$, then the
set of parameters where $\gamma$ is not transverse to the leaves of $\LL$
consists of isolated points.

\end{lemma}

This result immediately implies:

\begin{lemma} \label {aetra}

Let $f_\lambda$ be a non-trivial analytic family of unimodal maps.  Then
the set of non-regular Kupka-Smale parameters $\lambda_0$
such that $f_\lambda$ is not transverse to the topological class of
$f_{\lambda_0}$ at $\lambda_0$ is countable.

\end{lemma}

\subsection{Non-recurrent parameters}

The following result is due to Duncan Sands \cite {S}, but we will provide
a quick proof based on holomorphic motions and Lemma \ref {9.2}.

\begin{lemma} \label {nonre}

Let $f_\lambda$ be a non-trivial analytic family of unimodal maps.  Then
almost every parameter is regular or has a recurrent critical point.

\end{lemma}

\begin{pf}

If this is not the case, there would exist $\epsilon>0$ and
a set $X$ of parameters $\lambda$ of positive measure such that for $\lambda
\in X$,\\
(1)\, $\inf_{m \geq 1} |f_\lambda^m(0)|>\epsilon$ (by hypothesis),\\
(2)\, $f_\lambda$ is non-regular, Kupka-Smale and the critical orbit is
infinite (Lemma \ref {ks}).

Let us fix a density point $\lambda_0 \in X$ of $X$.  Using Lemma \ref
{5.1}, consider a nice interval $T=T[\lambda_0]=[-p,p] \subset
(-\epsilon,\epsilon)$ for $f_{\lambda_0}$, with $p$ preperiodic.
Let $T[\lambda]$, $\lambda-\lambda_0 \in (-\delta,\delta)$, $\delta>0$ small
denote the continuation of $T$.  Let $K[\lambda]$,
$\lambda-\lambda_0 \in (-\delta,\delta)$
denote the set of points in $I \setminus
T[\lambda]$ which never enter $T[\lambda]$ and do not
belong to the basin of hyperbolic attractors.

Since $K=K[\lambda_0]$ is an expanding set by Lemma \ref {5.2},
it persists in a complex neighborhood of $\lambda_0$:
there exists a family of homeomorphisms $h_\lambda:K \to \C$,
$\lambda \in \D_{\delta'}(\lambda_0)$, $\delta'<\delta$
depending continuously on $\lambda$, such that
$h_{\lambda_0}=\id$ and
$f_\lambda \circ h_\lambda=h_\lambda \circ
f_{\lambda_0}$.  It is easy to see (using Lemma \ref {5.2})
that for $\lambda \in \R$, $h_\lambda(K)=K[\lambda]$.
For each preperiodic orbit $p$ of $f$ in $K$,
it is clear that $\lambda \mapsto h_\lambda(p)$ is holomorphic in
$\D_{\delta'}(\lambda_0)$.
Since preperiodic orbits are dense in
$K$, it follows that $h[\lambda_0,\lambda] \equiv h_\lambda$
are actually transition maps of
a holomorphic motion $h$ over $\D_{\delta'}(\lambda_0)$.

Since $f_\lambda$ is non-trivial, $f_\lambda(0)$ does not belong to
$K[\lambda]$ for a dense set of $\lambda \in (-\delta,\delta)$,
so by Lemma \ref {9.2}, the path $\lambda \mapsto (\lambda,f_\lambda(0))$
is transverse to the leaves of $h$ outside of countably many
parameters $\lambda$.  So there exist parameters $\lambda \in X$
arbitrarily close to $\lambda_0$ which are density points of $X$ and
transversality points of the above path.  In order to avoid cumbersome
notation, let us assume that $\lambda_0$ is itself a transversality point.

It follows that there exists a real-symmetric quasiconformal
map $\chi$ (phase-parameter holonomy map\footnote{More precisely, $\chi$ is
obtained by applying first the local holonomy map between
the two transverse holomorphic curves
$\{\lambda_0\} \times V$ ($V$ a small neighborhood of
$f_{\lambda_0}(0)$)
and $\{(\lambda,f_\lambda(0))|\lambda \in
\D_{\delta'}(\lambda_0)\}$ followed by the projection in the first
coordinate.}) taking a neighborhood $V$
of $f_{\lambda_0}(0)$ to a neighborhood of $\lambda_0$, and
taking points in $K \cap V$ to parameters
$\lambda \in \chi(V)$ with $f_\lambda(0) \in K[\lambda]$.
In particular, $\chi(K \cap V) \supset X \cap \chi(V)$.

Since $K$ is an expanding set, it follows that there exists
$\rho>0$ such that in every $r$ neighborhood of $f_{\lambda_0}(0)$
there exists an interval of size at least $\rho r$ disjoint from
$K$.  Since $\chi|V \cap \R$ is quasisymmetric, this property is
preserved: there exists $\rho'>0$ such that in every $r$ neighborhood
of $\lambda_0$ there exists an interval of size at least $\rho' r$ not
intersecting $X$.  This contradicts the hypothesis that $\lambda_0$ is a
density point of $X$\footnote
{It is easy to see that this argument gives much more information on the
size of $X$.  One can see for instance that the Hausdorff
dimension of $X$ in $\lambda_0$ (defined as the infimum of the Hausdorff
dimension of $X \cap \D_\epsilon(\lambda_0)$)
is no greater than the Hausdorff dimension of
$K$ in $f_{\lambda_0}(0)$, which is known to be
less than $1$.  Notice that $X$ is essentially the set of non-regular
non-recurrent parameters avoiding a definite neighborhood of $0$.
We should remark that these ideas show also that the
Hausdorff dimension of the set of non-regular
non-recurrent parameters is usually
$1$ except for some trivial situations.}.
\end{pf}

\subsection{Infinitely renormalizable maps}

\begin{lemma} \label {infrenor}

Let $f_\lambda$ be a non-trivial analytic family of unimodal maps.  Then
the set of infinitely renormalizable parameters has Lebesgue measure zero.

\end{lemma}

\begin{pf}

Let $X$ be the set of parameters $\lambda$ such that $f_\lambda$ is
infinitely renormalizable, and let
$\lambda_0 \in X$ be a density point of $X$.
By Lemma \ref {5.6}, there exists a nice interval $T[\lambda]$,
$|\lambda-\lambda_0|<\delta$, which is periodic (of period, say, $m$)
such that $f^m|T[\lambda]$ has negative Schwarzian
derivative.  In particular, if $A_\lambda:T[\lambda] \to I$ is affine,
$g_\lambda=A_\lambda \circ f_\lambda^m \circ A_\lambda^{-1}$,
$|\lambda-\lambda_0|<\delta'$ is an analytic family of quasiquadratic maps,
which is non-trivial (because $f_\lambda$ is).
By Theorem B of \cite {ALM}, for almost every $\lambda$, $g_\lambda$ is not
infinitely renormalizable.  It is clear that if $\lambda \in X$ and
$|\lambda-\lambda_0|<\delta'$ then $g_\lambda$ is infinitely renormalizable,
so $\lambda_0$ is not a density point of $X$, contradiction.
\end{pf}

\subsection{Simple maps}

The following argument is adapted from the corresponding result of Lyubich
for the quadratic family \cite {parapuzzle}.

\begin{lemma} \label {simple}

Let $f_\lambda$ be a non-trivial analytic family of unimodal maps.  Then
almost every parameter $\lambda$ with $f_\lambda \in \FF$ is simple.

\end{lemma}

\begin{pf}

If this is not the case, we could find $C>0$, $\rho<1$,
$m \geq 0$ and
a set $X$ of parameters of positive
measure such that if $\lambda_0 \in X$ then

\begin{enumerate}

\item $f_{\lambda_0} \in \FF$ and is not simple (by hypothesis),

\item $f_\lambda$ is transverse at $\lambda_0$ (by Theorem \ref {aetra}),

\item The sequence of parameter windows $J_n[\lambda_0]$ associated to
$\lambda_0$ are defined for $n \geq m$
(by Theorem \ref {phpatra}),

\item If $n_{k,\lambda_0}-1$ denotes the sequence of non-central
levels of the principal nest of $f_{\lambda_0}$ then for
$n_{k,\lambda_0} \geq m$,
$|J_{n_{k,\lambda_0}+1}[\lambda_0]|/|J_{n_{k,\lambda_0}}[\lambda_0]|<C
\rho^k$ (by Theorem \ref {phpatra}).

\end{enumerate}

Consider now the set $X_k$, $k \geq m$
of parameters $\lambda_0 \in X$ such that the
return of level $n_{k,\lambda_0}$ is central.  Let $\Delta_k$
be the union of $J_{n_{k,\lambda_0}}[\lambda_0]$, $\lambda_0 \in X_k$
and $\Pi_k$ be the union of $J_{n_{k,\lambda_0}+1}[\lambda_0]$, $\lambda_0
\in X_k$.

Then each connected component $J_{n_{k,\lambda_0}}[\lambda_0]$
of $\Delta_k$ contains a single connected component
$J_{n_{k,\lambda_0}+1}[\lambda_0]$ of $\Pi_k$, and thus $|\Pi_k|/|\Delta_k|<C
\rho^k$, so that $|X_k| \leq |\Pi_k|<C \rho^k |\Delta_k| \leq C \rho^k
|\Delta_m|$.
On the other hand, $X \subset \cap_{k_0 \geq m} \cup_{k \geq k_0} X_k$ and
thus, $|X| \leq \inf_{k_0 \geq m} \sum_{k \geq k_0} C \rho^k |\Delta_m|=0$,
contradiction.
\end{pf}

The proof of Theorem A is concluded.

\section{Proof of Theorem B} \label {thmB}

We will give now a proof of Theorem B using a parameter exclusion argument. 
In the first version of this work (in \cite {Av}), a different proof
was given relying
on the refined statistical analysis of \cite {AM}, but we will
give a much simpler argument
based on a modified version of the quasisymmetric capacities of \cite {AM},
which allows us to get rid of distortion estimates and at the same time to
work with a fixed quasisymmetric constant.

\subsection{Measure estimate}

Define the modified $\g$-qs capacity of a set $X$ in an interval $I$ as
$$
p_\g(X|I)=\sup \frac {|h_1 \circ h_2(X \cap I)|} {|h_1 \circ h_2(I)|}
$$
where $h_1:\R \to \R$ is $\g$-qs and
$h_2:I \to \R$ is a diffeomorphism (onto its image)
with non-negative Schwarzian derivative.

Notice that if $F:T_1 \to T_2$ is a diffeomorphism with
non-positive Schwarzian derivative and $X \subset T_1$ then
$$
p_\g(X|T_1) \leq p_\g(F(X)|T_2).
$$
This is the main advantage of modified quasisymmetric capacities over the
traditional ones of \cite {AM}.

By the Koebe Principle, if $h:I \to I$ is a diffeomorphism and
has non-positive Schwarzian derivative then
$h([-\epsilon,\epsilon])=O(\epsilon)$.
By H\"older continuity of $\g$-qs maps, we get
$$
p_\g([-\epsilon,\epsilon]|[-1,1])=O(\epsilon^\kappa)
$$
for some $0<\kappa<1$ depending on $\g$.


For a map $f \in \FF$ with principal nest $\{I_n\}$, let $s$ be as in Lemma
\ref {5.7}, and let
$$
\alpha_n=p_\g(s(\cup I^j_n)|s(I_n)).
$$

Let us consider the components $T^k_n$ of $(R_{n-1}|I^0_{n-1})^{-1}
(\cup I^j_{n-1})$.
We reserve the index $0$ for the component containing $0$, and the indexes
$-1$ and $1$ for the components of $(R_{n-1}|I^0_{n-1})^{-1}(I^0_{n-1})$.
If $|k|>1$ then $R_{n-1}|T^k_n$ is a diffeomorphism onto some
$I^j_{n-1}$, $j \neq 0$ and $R_{n-1}^2|T^k_n$ is a diffeomorphism onto
$I_{n-1}$.  Let
$$
\epsilon_n=p_\g(s(\cup_{|k|>1} T^k_n)|s(I_n)).
$$

\begin{lemma} \label {10.1}

If $n$ is sufficiently large,
$(1-\alpha_{n+1}) \geq (1-\epsilon_{n+1})(1-\alpha_n)$.

\end{lemma}

\begin{pf}

If $|k|>1$ then $s(T^k_{n+1})$ is taken to
$s(I_n)$ by $s \circ R_n^2
\circ s^{-1}$ which has negative Schwarzian derivative for $n$ big.
In particular
$$
p_\g(s(\cup I^j_{n+1})|s(T^k_{n+1})) \leq p_\g(s(\cup C^\d_n)|s(I_n))
\leq \alpha_n.
$$
Thus $p_\g(s(\cup I^j_{n+1})|I_{n+1}) \leq
\epsilon_{n+1}+(1-\epsilon_{n+1}) \alpha_n$.
\end{pf}

\begin{lemma} \label {10.2}

If $f$ is simple then the $\epsilon_n$ decay exponentially fast.

\end{lemma}

\begin{pf}

If $f$ is simple then $|s(I_{n+1})|/|s(I_n)|$ decays exponentially fast
by Lemma \ref {5.5}.  In particular,
by the Koebe Principle, for each $j$, each of the
connected components of $s(I_{n+1} \setminus I^j_{n+1})$ is exponentially
(in $n$) bigger than $s(I^j_{n+1})$.  This implies that, for each $k$,
each component of $s(I_{n+2} \setminus T^k_{n+2})$ is exponentially bigger
than $s(T^k_{n+2})$ (using the Koebe Principle),
so $p_\g(s(T^k_{n+2})|s(I_{n+2}))$ decays exponentially
and so does $\epsilon_n$.
\end{pf}

\begin{lemma} \label {10.4}

If $f \in \FF$ does not admit a quasiquadratic renormalization then
$\cup I^j_n$ is not dense in $I_n$, for $n$ sufficiently big.

\end{lemma}

\begin{pf}

Up to considering a renormalization or unimodal restriction,
we may assume that $f$ is non-renormalizable and does not admit unimodal
restriction.  It is easy to see that if $x \in I$ never enters $I_1$ then
the iterates of $x$ accumulate on an orientation preserving
fixed point of $f$, and since $f$
does not admit a unimodal restriction, we conclude that $x \in \partial I$.

Since $f$ is not conjugate to a quadratic map,
there exists an interval $T$ whose orbit does not
accumulate on the critical point (Lemma \ref {quad}).
Let $n$ be biggest with the orbit of $T$ intersecting $I_n$ ($T$ intersects
$I_1$ by the previous discussion).
Of course, the set of points which land on $I_{n+1}$ does not intersect the
orbit of $T$, and so is not dense in $I_n$. 

It is easy to see that if the set of points in $I_m$ which eventually land
in $I_{m+1}$ is not dense in $I_n$ then $\cup I^j_{m+1}$ is not dense on
$I_{m+1}$.
In particular, by induction, $\cup I^j_m$ is not dense in $I_m$ for
$m \geq n+1$.
\end{pf}

\begin{lemma}

If $f$ does not admit a quasiquadratic renormalization then for $n$ large
enough, $\alpha_n<1$.

\end{lemma}

\begin{pf}

Let $n$ be large enough such that there exists an open interval
$E \subset I_n$ disjoint from $\cup I^j_n$, and
$s \circ R_n \circ s^{-1}$ has negative Schwarzian derivative.  We may
assume that $E \subset T$, where $\overline T \subset \inter I_n$
is a symmetric interval containing $I^0_n$.  By the Koebe Principle, there
exists $C>0$ such that if $h_2:s(I_n) \to \R$ has non-positive
Schwarzian derivative then $|h_2(s(E))|>C |h_2(s(T))|$.
In particular, there
exists $\epsilon>0$ such that if $h_1:\R \to \R$ is $\g$-qs, then, with
$h=h_1 \circ h_2$, we have $|h(s(E))|>\epsilon
|h(s(T))| \geq \epsilon |h(s(I^0_n))|$.

For $\d \in \Omega$, let
$E^\d=(R_n^\d)^{-1}(E)$.  Since $(R^\d_n)^{-1}$ has non-positive Schwarzian
derivative, we see that for any $h$ as above,
$|h(s(E^\d))|>\epsilon |h(s(C^\d_n))|$.
Notice that all the intervals $E^\d$, $\d \in \Omega$ are disjoint, and
$\cup E^\d$ does not intersect $\cup C^\d_n$ so
$$
p_\g(s(\cup C^\d_n)|s(I_n)) \leq \frac {1} {1+\epsilon}.
$$

\comm{
If $h_1:\R \to \R$ is a
$\g$-qs map and $h_2:s(I_n) \to \R$ has non-negative Schwarzian
derivative, we see that $h_2|s(I^j_n)$ has bounded distortion by the Koebe
Principle.  Notice that this implies that for some $\epsilon>0$,
if $j \neq 0$, $|h_1 \circ h_2(s(E^j))|>
\epsilon |h_1 \circ h_2 (I^j_n)|$.  Thus
$$
p_\g(s(I_n \setminus \cup E^j)|s(I_n)) \leq
p_\g(s(I^0_n)|s(I_n))+(1-p_\g(s(I^0_n)|s(I_n))(1-\epsilon)<1.
$$
}

By a previous argument of the proof of Lemma \ref {10.1},
$$
p_\g(s(\cup I^j_{n+1})|s(T^k_{n+1})) \leq p_\g(s(\cup C^\d_n)|s(I_n))<1
$$
for $|k|>1$.

Thus, $p_\g(s(\cup I^j_{n+1})|s(I_{n+1})) \leq
\epsilon_{n+1}+(1-\epsilon_{n+1})p_\g(s(\cup C^\d_n)|s(I_n))<1$.
\end{pf}

\begin{lemma} \label {densityhyp}

Let $f_\lambda$ be a one-parameter
non-trivial analytic family of unimodal maps
satisfying the Phase-Parameter relation at a parameter
$\lambda_0$ (in particular, $f=f_{\lambda_0}$ is simple).  Assume that $f$
does not admit quasiquadratic renormalization.
Then $\lambda_0$ is not a density
point of non-hyperbolic parameters\footnote{One can actually use those
techniques to show that $\lambda_0$ is a
density point of hyperbolic parameters, see Remark \ref
{compest} for the complex counterpart.}.

\end{lemma}

\begin{pf}

Let $J_n$ and $\Xi_n$ be as in the Topological Phase-Parameter relation.
Since $|J_n| \to 0$, and $\lambda_0 \in \Xi_n(I^{\tau_n}_n) \subset J_n$,
it is enough to show that
then there exists $\alpha<1$ such that
$$
\limsup \frac {|\Xi_n(\cup C^\d_n \cap I^{\tau_n}_n)|}
{|\Xi_n(I^{\tau_n}_n)|} \leq \alpha<1.
$$
Indeed, if $\lambda \notin \Xi_n(\cup
C^\d_n)$ then the critical point is non-recurrent.  By Lemma \ref {nonre},
for almost every non-recurrent parameter, $f_\lambda$ is hyperbolic.

Fix $1<\hat \g<\g$
By PhPa1, $\Xi_n|K_n \cap I^{\tau_n}_n$ is $\hat \g$-qs
for $n$ big enough.   On the other hand, for $n$ big enough,
$s^{-1}|s(I^{\tau_n}_n)$ is $C^1$ close to being linear
(because $s$ is analytic, and in particular $C^1$, and $s(I^{\tau_n}_n)$
is small).  So $\Xi_n \circ s^{-1}|s(K_n \cap I^{\tau_n}_n)$ is
$\g$-qs for $n$ big enough.  In particular
$$
\frac {|\Xi_n(\cup C^\d_n \cap I^{\tau_n}_n)|}
{|\Xi_n(I^{\tau_n}_n)|} \leq
\frac {|\Xi_n \circ s^{-1} s(\cup C^\d_n \cap I^{\tau_n}_n)|}
{|\Xi_n \circ s^{-1} s(I^{\tau_n}_n)|} \leq p_\g(s(\cup
C^\d_n)|s(I^{\tau_n}_n) \leq \alpha_n.
$$
By Lemmas \ref {10.1}, \ref {10.2} and \ref {10.4},
$\alpha=\limsup \alpha_n<1$.
\end{pf}

Theorem A and Lemma \ref {densityhyp} imply Theorem B for one-parameter
families.

\subsubsection{Many parameters} \label {many}

The argument of Lemma \ref {ks} implies the following:

\begin{lemma} \label {ksmany}

Let $\{f_\lambda\}_{\lambda \in \Lambda}$ be a $k$-parameter
non-trivial analytic family of unimodal maps.  The set of parameters which
are not Kupka-Smale or have a periodic or preperiodic critical point
is contained in a countable union of analytic submanifolds of $\Lambda$,
of codimension at least $1$, and so has Lebesgue measure zero.

\end{lemma}

Let us now show how the one-dimensional version of Theorem B implies the
general case.  Let $\{f_\lambda\}_{\lambda \in \Lambda}$
be a $k$-parameter analytic family of unimodal maps.  By Lemma \ref
{ksmany}, we just have to show that for any Kupka-Smale parameter
$\lambda_0 \in \inter \Lambda$, there exists a small
$\epsilon>0$, such that, letting $B_\epsilon \subset \Lambda$ be the ball
around $\lambda_0$ of radius $\epsilon$, almost
every parameter in $B_\epsilon$ is either regular or admits a quasiquadratic
renormalization.

Using Theorem \ref {8.1},
if $\epsilon$ is sufficiently small, $\lambda \mapsto f_\lambda$ is an
analytic map from $B_\epsilon$ to some open set $\VV$ where the hybrid
lamination $\LL$ is defined.  Let $\lambda_1 \in B_\epsilon$ be a regular
parameter.  If $L$ is a line in $\R^k$ through $\lambda_1$, then by Lemma
\ref {9.2}, $L \cap B_\epsilon$ is not contained in the topological class
of a non-regular parameter, and so regular parameters are dense in $L \cap
B_\epsilon$.  By the one-dimensional Theorem B, we see that almost every
non-regular parameter in $L \cap B_\epsilon$ is quasiquadratic.
By Fubini's Theorem, almost every non-regular
parameter in $B_\epsilon$ is quasiquadratic.

This completes the proof of Theorem B.

\section{Proof of corollaries} \label {pf}

\subsection{Some conditions related to good ergodic properties}

Let us first recall the conditions on the critical orbit stated in the
introduction.  Let $f \in \U^2$.  We say that $f$ is {\it Collet-Eckmann}
if the lower Lyapunov exponent of the critical value is bigger than zero:
\begin{equation}
\liminf \frac {\ln |Df^n(f(0))|} {n}>0.
\end{equation}
We say that $f$ has {\it subexponential recurrence} if
\begin{equation}
\limsup \frac {-\ln |f^n(0)|} {n}=0.
\end{equation}
We say that $f$ has {\it polynomial recurrence} if
\begin{equation}
\g=\limsup \frac {-\ln |f^n(0)|} {\ln (n)}<\infty,
\end{equation}
and in this case, we call $\g$ the {\it exponent} of the recurrence.

We introduce the following additional condition: $f$ is called
{\it Weakly Regular} if
\begin{equation}
\lim_{\delta \to 0} \liminf_{n \to \infty} \frac {1} {n}     
\sum_{\ntop {1 \leq k \leq n} {f^k(0) \in (-\delta,\delta)}}
\ln |Df(f^k(0))|=0.
\end{equation}

Notice that polynomial recurrence is much stronger than subexponential
recurrence.

\begin{rem} \label {ergodic}

Maps satisfying the Collet-Eckmann and the subexponential recurrence
conditions have been intensively studied after the works of
Benedicks and Carleson.  Those two conditions give a very
precise control of the critical orbit.  They are not sufficient to show that
$f$ has good statistical properties however: one must also ask that $f$ has
a renormalization with all periodic orbits repelling (which is then
conjugate to
a quadratic polynomial).  Under this additional
assumption, it is possible to show that $f$ has an absolutely continuous
invariant measure (see \cite {BY}).

In order to study further
the properties of $\mu$, it is convenient to consider the smallest periodic
nice interval $T$ of $f$ ($f$ is not infinitely renormalizable,
since it has an absolutely continuous invariant measure).
The first return map $f^m:T \to
T$ can be then rescaled to a unimodal map $\hat f$, which also possess an
absolutely continuous invariant measure $\hat \mu$.

Assuming that $f$ is also Kupka-Smale and using Lemma \ref {5.2}, we see
that the dynamics of $f$ splits in a hyperbolic part, that describes
points $x \in I$ which never enter $\inter T$, and an interesting part
described by $\hat f$.

The measurable dynamics of $\hat f$ are described by
$\hat \mu$: for almost every $x \in I$ and any continuous function $\phi:I
\to \R$ we have
$$
\frac {1} {n} \sum_{k=0}^{n-1} \phi(\hat f^k(x))=\int \phi d\hat \mu.
$$

Since $\hat f$ is non-renormalizable, it follows that $\hat \mu$ is
supported on $[\hat f^2(0),\hat f(0)]$, and $(\hat f,\hat \mu)$ is
exponentially mixing\footnote{For a certain class of observables, for
instance, of bounded variation.} (see \cite {Y2}).

The condition of Weak Regularity is important to show that $(\hat f,
\hat \mu)$ is stochastically stable\footnote{For a certain class of
i.i.d. absolutely continuous
stochastic perturbations, the perturbed system possess a stationary
measure which is close to $\hat \mu$ in the weak topology.}
(see \cite {T2}).  If we assume a little bit
more smoothness, $f \in \U^3$, the Weak Regularity condition
is not necessary, and it is possible to show that
$(\hat f,\hat \mu)$ is stochastically stable in a stronger
sense\footnote{Densities of stationary measures of perturbed
systems are close to the density of $\hat \mu$ in the $L^1$ sense.}
(see \cite {BV}).

\end{rem}

\subsection{Analytic families}

We will actually prove the following result, which is a more precise form of
Corollaries C and E:

\begin{thm} \label {C precise}

Let $f_\lambda$, be a non-trivial analytic
family of unimodal maps.  Then almost every
non-regular parameter is Kupka-Smale and
has a quasiquadratic renormalization which satisfies the Collet-Eckmann
condition and is polynomially recurrent with exponent 1.

\end{thm}

\begin{pf}

We will prove the stated result for one-parameter families, the general case
reducing to this one by the argument of \S \ref {many}.

By Theorems A and B of \cite {AM}, the conclusion of the theorem holds
for the quadratic family.  However, the only properties of the quadratic
family that are actually used in the proof is that it is an analytic family
of quasiquadratic maps with negative Schwarzian derivative for which
the Phase-Parameter relation holds at almost every parameter, see Remark 3.3
of that paper.  Due to the work of Kozlovski, the hypothesis of negative
Schwarzian derivative can also be removed (this can be checked directly
using Lemma \ref {5.7}).  Using our Theorem A, we get the result
for analytic families of quasiquadratic maps.

\comm{
We already know that for almost every non-regular parameter $\lambda_0$,
$f=f_{\lambda_0}$ is simple and satisfies the phase-parameter relation
(by Theorem A), and has a renormalization which is
quasiquadratic (Theorem B).  By Lemma \ref {aetra} we also know that
$f_\lambda$ is transverse to the topological class of $f$ at $\lambda_0$.
}

Let us now consider the general case.  By Theorem A, almost every
non-regular parameter is simple, and by Theorem B, almost every non-regular
parameter has a quasiquadratic renormalization.
Let us fix such a parameter $\lambda_0$.

Let $T$ be the smallest periodic nice interval for $f_{\lambda_0}$
(of period $m$).  For $\lambda$ near $\lambda_0$, the interval
$T$ has a continuation $T[\lambda]$.  Consider
the analytic family $\hat f_\lambda=A[\lambda] \circ f_\lambda^m \circ
A[\lambda]^{-1}$, $|\lambda-\lambda_0|<\epsilon$,
where $A[\lambda]:T[\lambda] \to I$ is affine.  Then $\hat f_\lambda$ is
$C^\infty$ close to $\hat f_{\lambda_0}$, which is quasiquadratic,
so we conclude that for $\epsilon>0$ small, $\hat f_\lambda$,
$|\lambda-\lambda_0|<\epsilon$
is an analytic family of quasiquadratic maps.
Since $f_\lambda$ is non-trivial, $\hat f_\lambda$ is also non-trivial.

In particular, by the quasiquadratic case,
for almost every $\lambda$ near $\lambda_0$, $\hat f_\lambda$
is either regular or satisfy
the Collet-Eckmann condition and its critical point is polynomially
recurrent with exponent $1$.  In particular, the same holds for $f_\lambda$,
which concludes the proof of the theorem.
\comm{
It is clear that$f_\lambda$ is
non-trivial, so is $\hat f_\lambda$.  We conclude that almost every parameter
$\lambdanon-regular
parameter is an analytic family of quasiquadratic maps.  It also satisfies the
phase-parameter relation at almost every non-regular parameter.  This is
enough to repeat the proof of Theorems A and B of \cite {AM}.
}
\end{pf}

\begin{rem}

Notice that the proof of Theorem A in \cite {AM2} could not use directly the
proof of \cite {AM} (the argument needs modifications which are dealt in
the Appendix of \cite {AM2}), since their main phase-parameter tool
essentially amounts to comparing
the phase-space of a non-trivial family with the parameter space of the
quadratic family.  This distorts the estimates and makes it impossible to
obtain the exponent of the recurrence.

\end{rem}

\subsection{Smooth families}

Recall that if $\Lambda \in \R^k$ is a bounded open connected
domain with smooth boundary,
$\UF^r(\Lambda)$ is the space of $C^r$ families of unimodal maps
parametrized by $\Lambda$, and is a Baire space.

\begin{thm} \label {D}

Let $f_\lambda$, $\lambda \in \Lambda$
be a non-trivial family of unimodal maps.  For every
$\epsilon>0$ there exists a neighborhood $\VV \subset \UF^2(\Lambda)$ of
$f_\lambda$ such that if $g_\lambda \in \VV$
then, outside a set of parameters $\lambda$
of measure at most $\epsilon$, $g_\lambda$ is either
regular or is Kupka-Smale and has a renormalization
with all periodic orbits repelling satisfying the
Collet-Eckmann, subexponential recurrence, and Weak Regularity conditions.

\end{thm}

\begin{pf}

Using Vitali's Covering Lemma, let $\{B_i\}$, $\{C_i\}$ be
finite families of disjoint
closed balls covering the parameter space up to a set of Lebesgue
measure $\epsilon/2$ such that:\\
(1)\, For $\lambda \in B_i$, $f_\lambda$ is regular;\\
(2)\, For $\lambda \in C_i$, there exists a nice
interval $T_i[\lambda]$, which is periodic of period $m_i$,
depending continuously on
$\lambda$ such that
$f_\lambda^{m_i}:T_i[\lambda] \to T_i[\lambda]$ can be rescaled
to a quasiquadratic map $\hat f_{i,\lambda}$.

It is easy to see that if $g_\lambda$ is $C^2$ close to $f_\lambda$,
then:\\
(1)\, For every $\lambda \in B_i$, $g_\lambda$ is regular;\\
(2)\, For every $\lambda \in C_i$, there exists an
interval $T^g_i[\lambda]$, depending continuously on $\lambda$, close to
$T^g_i[\lambda]$, such that
$g^{m_i}_\lambda:T^g_i[\lambda] \to T^g_i[\lambda]$ can be rescaled to a
unimodal map $\hat g_{i,\lambda}$, and the family $\hat g_{i,\lambda}$ is
$C^2$ close to $\hat f_{i,\lambda}$.

The family $\hat f_{i,\lambda}$ is non-trivial,
so by Theorem B of \cite {ALM}, the set of parameters in $C_i$ such that
$\hat g_{i,\lambda}$ is either regular or has all periodic orbits repelling
and satisfies the Collet-Eckmann, subexponential recurrence, and
Weak Regularity conditions, has Lebesgue measure at least
$|C_i|(1-\epsilon/4)$, provided $g_\lambda$ is close enough to
$f_\lambda$.  The result follows.
\end{pf}

\begin{rem}

In particular, if $f_\lambda$ is a non-trivial analytic family of unimodal
maps, almost every parameter is Weakly Regular.

\end{rem}

Recall that by Remark \ref {dense},
non-trivial analytic families are dense in
$\UF^r(\Lambda)$.  Using Theorem \ref {D} and an easy Baire argument we
obtain the following precise version of Corollary D:

\begin{thm} \label {D precise}

In a generic family $f_\lambda$ in $\UF^r(\Lambda)$, $r=2,...,\infty$ for
almost every non-regular
parameter $\lambda_0 \in \Lambda$, $f=f_{\lambda_0}$
is Kupka-Smale and has a renormalization which has all periodic
orbits repelling and satisfies the Collet-Eckmann, subexponential
recurrence, and Weak Regularity conditions.

\end{thm}

\appendix

\section{Hybrid classes}

In this section we will give a global characterization of the leaves of the
lamination $\LL$ of Theorem \ref {8.1}.

Notice that the leaves of $\LL$ are claimed
to coincide with topological classes only in the non-regular case:
the partition in topological classes is not a lamination because
regular topological classes are open sets.  It turns out that the behavior
of the regular leaves of $\LL$ can be quite arbitrary.  In order to give a
global characterization of the leaves of $\LL$, we need
to introduce once and for all an arbitrary, but fixed, way
to refine the topological classes of regular maps.  We shall call this
refinement the {\it hybrid lamination}.

If $f$ is non-regular, the hybrid class of $f$ is just the set of all
non-regular maps $g$ which are topologically conjugate to $f$.

Let $f$ be a regular map, and let $A$ be the set of attracting periodic
orbits of $f$ and let $B=\{x \in I|f^n(x) \to A\}$
denote the basins of the attracting periodic orbits of $f$.
Notice that if $f$ is a regular map, there exists
a minimal $m \geq 0$ such that $f^m(0)$ belongs to a periodic connected
component of $B$.  It is possible to show that if $f$ is quasiquadratic,
then $m=0$.  It turns out that if $m=0$ (this case will be
called essential), there is a natural way to refine the topological
class of $f$: the hybrid class of $f$ is the set of all regular maps
$g$ which are topologically conjugate to $f$ and the multiplier
of the periodic orbit that attracts $0$ is the same for both maps (this
definition agrees with the one of \cite {ALM} in the quasiquadratic case).

In the non-essential case, there is no natural way to refine the topological
class of $f$, so we fix an arbitrary way that works.

\begin{definition}

Let $f$ be a Kupka-Smale map.  We say that a homeomorphism
$h:I \to \C$ is $f$-admissible if the following holds.
Let $T$ be a periodic component of $B \setminus A$ which does not contain
$0$, and, writing $T=(a,b)$ with $|a|<|b|$, we have that the interval
$[-a,a]$ is nice.  Then $h$ takes $d=(a+b)/2$ to $h(d)=(h(a)+h(b))/2$ and
$h|[d,f^q(d)]$ is affine, where $q$ is the period of $T$.

\end{definition}

\begin{definition}

Let $f$ be a regular map of non-essential type.
The hybrid class of $f$ is defined as the
set of all regular maps $g$ such that there exists an $f$-admissible
topological conjugacy between $f$ and $g$.

\end{definition}

The following proposition
is elementary, and shows that the definition of hybrid
class is minimally adequate:

\begin{prop}

Let $f$ be a regular map.  Then its hybrid class intersects $\UU_a$ in a
codimension-one analytic submanifold.

\end{prop}

Moreover, with this definition, it is possible to prove the full
Theorem \ref {8.1} in the case of hyperbolic maps $f$.  The case of
infinitely renormalizable $f$ can be dealt by reduction to the quasiquadratic
case using renormalization (dealt in Theorem A of \cite {ALM}),
see Lemma \ref {5.6}.

\subsection{Persistent puzzle}

The remaining case of Theorem \ref {8.1} is trickier and one needs to go
into the proof of \cite {ALM}.  We will discuss here only the main
modification one needs to make in order to adapt the argument.  This
modification concerns the main tool used in the finitely renormalizable
case, the concept of persistent puzzle, whose definition needs to be
adapted.  We follow basically the approach of \cite {Av}.


Assume that $f \in \FF$.  As in \S \ref {6.1},
fix a level $\n$ of the principal nest and assume
that $|I_\n|/|I_{\n-1}|$ is very small.  Let us consider the
first landing map to $A^0=I_\n$, the connected components of its domain
are denoted $A^j$.  Let $A^1$ be the component of $f(0)$, and let
$A^1=[l,r]$, with $l<r$.  Let $V^j$ be the complexification of the $A^j$
obtained as in Lemma \ref {geometric puzzle}.  Let $V$ be the union
of all $V^j$ such that $V^j \cap \R \subset [-1,r]$.  We shall
informally call $V$ the {\it puzzle}.

Let $\VV \subset \AAA_a$ be a real-symmetric
neighborhood of $f$.  We will say that
the puzzle {\it persists} in $\VV$ if there exists a real-symmetric
holomorphic motion $h$ over $\VV$ given by a family of
transition maps $h[f,g]=h_g:\C \to \C$, $g \in \VV$ such that:

\begin{enumerate}

\item $h_g|\C \setminus \Omega_a=\id$;

\item $g \circ h_g|V \setminus V^0=h_g \circ f$,
$g \circ h_g|\partial V^0=f$;

\item $h_g|I$ is $f$-admissible
and $g \circ h_g|([-1,r] \setminus V)=h_g \circ f$.

\end{enumerate}

The following plays the role of Lemma 5.6 of \cite {ALM}.

\begin{lemma}

Let $f \in \FF \cap \UU_a$.
If $|I_\n|/|I_{\n-1}|$ is sufficiently small, then there exists a
neighborhood of $f$ where the puzzle persists.

\end{lemma}

The proof is the same as of Lemma 5.6 of \cite {ALM}, and we will not
reproduce the whole argument here, but only comment the main steps:

(1)\, One considers a holomorphic motion $h'$ of
$[-1,r] \setminus V$ which is $f$-admissible and equivariant: $g \circ
h'_g=h'_g \circ f$ (this
holomorphic motion exists because the dynamics of $f|[-1,r] \setminus V$
is hyperbolic) over a small neighborhood of $f$.

(2)\, Using the Canonical Extension Lemma, we extend $h'$ to a holomorphic
motion defined also on $\partial f(V_0)$.  Considering a slightly smaller
neighborhood $\VV'$ of $f$ we may extend $h'$ to $\C \setminus \Omega$ as
$\id$.

(3)\, One considers a holomorphic motion $h^0$ of $\overline V^0$
such that $g \circ h^0_g|\partial V^0=h'_g \circ f$ over a neighborhood
$\VV^0$ of $f$.

(4)\, One notices that for each $V^i$, $i \neq
0$, we can define (uniquely) a holomorphic motion $h^i$ on $V^i$
as a lift of $h^0|V^0$ over a small neighborhood $\VV^i$ of $f$.

(5)\, The (countably many) holomorphic motions $h'$, $h^i$ are
defined apriori over different neighborhoods of $f$,
but using again hyperbolicity of $f|[-1,r] \setminus V$, one sees that all
those holomorphic motions are defined over a definite neighborhood of $f$.

(6)\, An estimate of hyperbolic geometry shows that the several regions of
definition of those different holomorphic motions cannot collide in a
slightly smaller neighborhood of $f$, so they
define a common holomorphic motion which can be completed using the
Canonical Extension Lemma and satisfies automatically (1), (2), and (3).

\begin{rem}

The last condition of the definition of persistence defines
uniquely $h_g$ in $[-1,r] \setminus \overline V$.
This set is empty in the quasiquadratic case (and so this condition
does not appear in \cite {ALM}).  This (obvious) observation concerning
the first step is the only formal difference in the proof, the remaining
steps do not need to be modified.

\end{rem}

\begin{rem}

If $f$ is a Kupka-Smale, finitely-renormalizable, non-hyperbolic map, with a
non-recurrent critical point, a similar construction can be made.  In this
case, we take $T \subset T'$ nice intervals with preperiodic boundary such
that $0$ does not return to $T'$ and $|T|/|T'|$ is very small.  We let
$A^0=T$, and put $A^1$ as a domain of the first landing map to $A^0$ which is
contained in $[f(0),f(0)+\epsilon]$, $\epsilon$ very small.

\end{rem}

\begin{rem}

If $g_1,g_2 \in \VV \cap \UU_a$ are regular maps in the same hybrid class
then they are of non-essential type if and only if for all $m$ sufficiently
big,
$$
h_{g_1}^{-1}(g_1^m(0)),
h_{g_2}^{-1}(g_2^m(0)) \notin [-1,r] \setminus \overline V
$$
(use the Schwarz Lemma).  The definition of hybrid class implies
$$
h_{g_1}^{-1}(g_1^m(0))=h_{g_2}^{-1}(g_2^m(0)).
$$
This is important for
the application of the several pullback arguments of \cite {ALM}.

\end{rem}

One obtains Theorem \ref {8.1} in the finitely renormalizable, non-regular
case by repetition of the proof of
Theorem A of \cite {ALM}, taking into consideration the above remarks.

\section{Non-renormalizable parameters in the Mandelbrot set}

Let $p_c=z^2+c$ and let $\MM$ (the Mandelbrot set)
be the set of parameters $c \in \C$ such that
the orbit of $0$ does not escape to infinity under iteration by $p_c$.
The aim of this appendix is to show how the idea of the proof of
Theorem B can be coupled with Lyubich's result of \cite {parapuzzle}
to obtain the following theorem:

\begin{thm} \label {nren}

Let $\NR$ be the set of non-renormalizable quadratic
parameters with recurrent critical point and no indifferent
periodic orbits in the boundary of the Mandelbrot set.  Then $\NR$
has Lebesgue measure $0$.

\end{thm}

Theorem \ref {nren} implies easily Shishikura's
Theorem F stated in the introduction.

\comm{
\begin{thm} \label {infren}

The set of parameters in the boundary of the Mandelbrot set which are not
infinitely renormalizable has Lebesgue measure $0$.

\end{thm}
}

\begin{rem}

The reduction of Theorem F to Theorem \ref {nren} is obtained
using the following three steps:

(1)\, It is easy to pass from the non-renormalizable case to
the finitely renormalizable case
using renormalization techniques: the (countably many)
little copies of the Mandelbrot set are
related by renormalization to the original Mandelbrot set by a
quasiconformal (and thus absolutely continuous)
transformation, see \cite {universe}.
Alternatively, we can also repeat the
proofs for the little Mandelbrot copies.

(2)\, Quadratic polynomials with a neutral fixed point are contained in the
boundary of the main cardioid of the Mandelbrot set, which is a real
analytic curve (with one singularity) and thus has Lebesgue measure zero.

(3)\, The case of non-recurrent non-renormalizable polynomial
without neutral fixed points
can be treated easily using holomorphic motions, see our proof of
Lemma \ref {nonre} (it is enough to use that under those conditions
the set of points that never enter a small neighborhood of $0$
is a hyperbolic set
and thus persistent\footnote{This actually holds for any non-renormalizable
quadratic polynomial without neutral fixed points.}).

\end{rem}

To prove Theorem \ref {nren} we will make use of the Phase-Parameter
estimates described in Lemma~\ref {cphpa1}
and Lyubich's parapuzzle estimate (Theorem~\ref {p_lambda}).
Then, we will redo the estimates of Theorem B in the
complex setting to show that non-renormalizable parameters have Lebesgue
measure zero, because the critical point has a tendency to fall in the basin
of infinity (in the same way that in the real setting the critical
point has a tendency to fall in the basin of non-essential attractors).

\comm{
To prove Theorem \ref {nren} we will make use of the complex
phase-parameter relation of
\cite {AM} for the case of complex parameter values, which we will formulate
in the next section.  Then, in \S we will redo the estimates of \S in this
complex setting to show that non-renormalizable parameters have Lebesgue
measure $0$, because the critical point has a tendency to fall in the basin
of infinity (in the same way that in the real setting the critical
point has a tendency to fall in the basin of non-essential attractors).
}

\begin{rem}

Lyubich has another proof of Theorem \ref {nren}, also based on
\cite {parapuzzle} and estimates on the area of the set of points that
return to deep puzzle pieces.
Graczyk and Swiatek have also obtained a different proof
of Shishikura's Theorem.

\end{rem}

\subsection{Parapuzzle notation}

Let us fix $c_0 \in \NR$.  By Theorem \ref {p_lambda}, there exists a
neighborhood $\Lambda_1 \subset \C$ of $c_0$ and domains
$0 \in U_1[\lambda] \subset \C$, $\lambda \in \Lambda_1$ such that
the first return map to
$U_1[\lambda]$ by $p_\lambda$ induces a full $R$-family over
$\Lambda_1$.
To prove Theorem \ref {nren}, it is clearly sufficient to show that
$\Lambda_1 \cap \NR$ has Lebesgue measure zero.

For $\lambda \in \NR \cap \Lambda_1$, we can define a $R$-chain over
$\lambda$ since the critical point is recurrent.  Let us denote the
parameter domains of this chain by $\Lambda_i[\lambda]$.
Let $\NR^\infty \subset \NR \cap \Lambda_1$ be the set of
parameters $\lambda$
such that the chain $\RR_i$ over $\lambda$ has infinitely many central
levels, and let $\NR^0$ be the complementary set in $\NR \cap \Lambda_1$.
By Theorem \ref {3.5}, there exists a constant $C(\lambda)>0$, $\lambda \in
\NR \cap \Lambda_1$ such that $\mod(\Lambda_{n_k}[\lambda] \setminus
\overline {\Lambda_{n_k+1}[\lambda]}) >
C(\lambda) k$, where $n_k-1$ counts the
non-central levels of the chain.  If $\lambda \in \NR^0$, we actually have
linear growth of moduli (without passing through a subsequence), and by
Lemma~\ref {cphpa1}, conditions CPhPa1 and CPhPh1 are satisfied (with the
dilatation parameter $\g$ arbitrarily close to $1$) for $i$ sufficiently
big.

\subsection{Finitely many central cascades}

The argument of Lyubich which shows that almost every
real quadratic maps in $\FF$ is simple
applies in the complex setting and gives:

\begin{lemma}

$|\NR^\infty|=0$.

\end{lemma}

\begin{pf}

Let $\NR^\infty_\epsilon$ be the set of parameters $\lambda \in NR^\infty$
such that $C(\lambda) \geq \epsilon$. 
If $\NR^\infty$ has positive Lebesgue measure then we can select
$\epsilon$ such that $\NR^\infty_\epsilon$
also has positive Lebesgue measure.  Let $\NR^\infty_\epsilon(k)
\subset \NR^\infty_\epsilon$ be the set of parameters such that the
$n_k$ level is central.  If $\lambda \in \NR^\infty_\epsilon(k)$,
$\NR^\infty_\epsilon(k) \cap \Lambda_{n_k}[\lambda] \subset
\Lambda^0_{n_k}[\lambda]$,
thus $|\NR^\infty_\epsilon(k) \cap \Lambda_{n_k}[\lambda]|
\leq |\Lambda^0_{n_k}[\lambda]|$.

Since $C(\lambda) \geq \epsilon$, there exists $\delta$ and $k_0$ which only
depend on $\epsilon$ such that if $k>k_0$ then $\Lambda_{n_k}[\lambda]
\setminus \overline {\Lambda^0_{n_k}[\lambda]}$
contains a round annuli of moduli $k
\delta$.  This implies that
$|\Lambda^0_{n_k}[\lambda]| \leq e^{-k \delta'} |\Lambda_{n_k}[\lambda]|$
for some $\delta'$ depending on $\delta$.
For each $k$, the domains $\Lambda_{n_k}[\lambda]$, $\lambda \in
\NR^\infty_\epsilon(k)$ are either equal or disjoint,
and their union has Lebesgue measure at most $|\Lambda_1|$, so
$|\NR^\infty_\epsilon(k)|$ decays exponentially on $k$.
It follows immediately that
$\NR^\infty_\epsilon=\cap_{k \geq 1} \cup_{n \geq k} \NR^\infty_\epsilon(k)$
has Lebesgue measure zero, contradiction.
\end{pf}

\subsection{Area estimate}

Let $U$ be a bounded open set of $\C$ and $Z$ be a measurable set of $\C$.
Let
$$
c_\g(Z|U)=\sup \frac {|h(Z \cap U)|} {|h(U)|}
$$
where $h$ ranges over all quasiconformal homeomorphisms
$h:U \to \C$ with dilatation bounded by $\g$
and such that $h(U)$ is bounded.  The following two properties are
immediate:

\begin{enumerate}

\item If $V^j \subset U$ are disjoint open subsets
and $Z \subset \cup V^j$ then
$$
c_\g(Z|U) \leq \sup_j c_\g(Z|V^j) c_\g(\cup V^j|U).
$$

\item If $A,B \subset U$ are disjoint open subsets and
$Z \subset A \cup B$ then
$$
c_\g(Z|U) \leq c_\g(B|U)+(1-c_\g(B|U))c_\g(Z|B).
$$

\end{enumerate}

Denote by $V^k_n[\lambda]$
the connected components of the preimages of
$$
(R_{n-1}[\lambda]|U^0_n[\lambda])^{-1}(\cup U^j_{n-1}[\lambda]).
$$
We reserve the index $0$ for the component of $0$, so that $0 \in V^0_n$.
We also reserve the indexes $-1$ and
$1$ for the components of the preimages of $U_n[\lambda]$.

Fix some $\g>1$.  Let
\begin{align}
\epsilon_n(\lambda)&=c_\g(\cup_{|k| \leq 1} V^k_n[\lambda]|U_n[\lambda])\\
\alpha_n(\lambda)&=c_\g(\cup_j U^j_n[\lambda]|U_n[\lambda]).
\end{align}

\begin{lemma} \label {B.4}

Let $\lambda \in \NR^0$.  Then $\alpha_2<1$.

\end{lemma}

\begin{pf}

Notice that $\cup U^j_1[\lambda]$ is not dense in $U_1[\lambda]$ (otherwise
the filled-in Julia set of $p_\lambda$ would have to contain
$U_1[\lambda]$, but in our situation the filled-in Julia set of $p_\lambda$
has empty interior).  Thus, there exists a domain
$U^0_1[\lambda] \subset D[\lambda] \subset U_1[\lambda]$ such that
$U_1[\lambda] \setminus \overline
{D[\lambda]}$ is an annulus, and a non-empty open set $E[\lambda]
\subset D[\lambda] \setminus \cup U^j_1[\lambda]$.
By the Koebe distortion Lemma, if
$h:U_1[\lambda] \to \C$ is a $\g$-qc map with bounded image then
$|h(E[\lambda])|>C|h(U^0_1[\lambda])|$ for some constant $C>0$.

For $\d \in \Omega$, let
$E^\d[\lambda]=(R^\d_1[\lambda])^{-1}(E[\lambda])$.  We conclude that, for
any $\g$-qc map $h:U_1[\lambda] \to \C$ with bounded image,
we have $|h(\cup
E^\d[\lambda])|>C|h(\cup W_1^\d[\lambda])|$, so $c_\g(\cup
W_1^\d[\lambda]|U_1[\lambda])<1$.

If $|k|>1$ then $R_1^2[\lambda]|V^k_2[\lambda]$ is a diffeomorphism
onto $U_1[\lambda]$ and we conclude that
$c_\g(\cup U^j_2[\lambda]|V^k_2[\lambda])=c_\g(\cup
W_1^\d[\lambda]|U_1[\lambda])$.

Thus
$c_\g(\cup U^j_2[\lambda]|U_2[\lambda]) \leq \epsilon_2+(1-\epsilon_2)
c_\g(\cup W_1^\d[\lambda]|U_1[\lambda])<1$.
\end{pf}

\begin{lemma} \label {B.5}

If $\lambda \in \NR^0$ then $\epsilon_n(\lambda) \to 0$ exponentially fast.

\end{lemma}

\begin{pf}

Notice that if $R_{n-1}[\lambda](V^k_n[\lambda])=U^j_{n-1}[\lambda]$ then
$$
\mod(U_n[\lambda] \setminus \overline
{V^k_n[\lambda]}) \geq
\mod(U_{n-1}[\lambda] \setminus \overline
{U^j_{n-1}[\lambda]})/3,
$$
$$
\mod(U_{n-1}[\lambda] \setminus \overline
{U^j_{n-1}[\lambda]}) \geq \mod(U_{n-2}[\lambda] \setminus \overline
{U^0_{n-2}[\lambda]})/2.
$$
For $\lambda \in \NR^0$,
$\mod(U_{n-2}[\lambda] \setminus \overline
{U^0_{n-2}[\lambda]})$ grows linearly in $n$, so
$\inf_k \mod(U_n[\lambda] \setminus \overline
{V^k_n[\lambda]})$ also grows linearly, and this implies exponential decay
of $\sup_k c_\g(V^k_n[\lambda]|U_n[\lambda])$,
which implies exponential decay of $\epsilon_n$.
\end{pf}

\begin{lemma}

If $\lambda \in \NR^0$ then $\alpha(\lambda)=\sup_{n \geq 2}
\alpha_n(\lambda)<1$.

\end{lemma}

\begin{pf}

Indeed, if $|k|>1$ then $R^2_n[\lambda]|V^k_{n+1}[\lambda]$ is a
diffeomorphism onto $U_n[\lambda]$.  In particular,
$c_\g(\cup U^j_{n+1}[\lambda]|V^k_{n+1}[\lambda]) \leq c_\g(\cup
U^j_n[\lambda]|U_n[\lambda])=\alpha_n(\lambda)$.
Thus
$$
c_\g(\cup U^j_{n+1}[\lambda]|U_{n+1}[\lambda] \setminus \overline
{\cup_{|k| \leq 1} V^k_{n+1}[\lambda]}) \leq \alpha_n(\lambda),
$$
which implies
$\alpha_{n+1}(\lambda) \leq
\epsilon_{n+1}(\lambda)+(1-\epsilon_{n+1}(\lambda))\alpha_n(\lambda)$
and
$$
1-\alpha_{n+1}(\lambda) \geq
(1-\epsilon_{n+1}(\lambda))(1-\alpha_n(\lambda)).
$$

If $\lambda \in \NR^0$,
$\epsilon_n(\lambda)$ decays exponentially (Lemma \ref
{B.5}) and $\alpha_2(\lambda)<1$ (Lemma \ref {B.4}), so the result follows.
\end{pf}

If $\NR^0$ has positive measure,
there exists $\alpha>0$, $k>0$ and a positive measure set $X$
such that for $\lambda \in X$,
$\alpha(\lambda)<\alpha$ and for $n>k$
the estimate CPhPa1 of Lemma~\ref {cphpa1} is valid
with a constant smaller than $\g$.

Let $Y \supset X$ be an open set such that
$\alpha |Y|<|X|$.  For every parameter
$\lambda \in X$, let $\mu(\lambda)$
be the smallest $j>k$ such that
$\lambda \in Z[\lambda]=\Lambda^{\tau_j(\lambda)}_j[\lambda]
\subset Y$ (such a $j$ exists since $\cap \Lambda_j[\lambda]=\{\lambda\}$).
The resulting collection of parameter domains $Z[\lambda]$, $\lambda \in
X$ are either disjoint or equal. 
To reach a contradiction, it is enough to show that
$\alpha |Z[\lambda]| \geq
|X \cap Z[\lambda]|$,
for in this case $\alpha |Y| \geq |X|$.  But this is an immediate
consequence of CPhPa1, for
$$
\frac {|X \cap Z[\lambda]|} {|Z[\lambda]|} \leq
c_\g(\cup W^\d_{\mu(\lambda)}|
U^{\tau_{\mu(\lambda)}(\lambda)}_{\mu(\lambda)}) \leq
c_\g(\cup U^j_{\mu(\lambda)}|U_{\mu(\lambda)}) \leq \alpha,
$$
since $\tau_{\mu(\lambda)} \neq 0$ by hypothesis (notice that we even have
$|\MM \cap Z[\lambda]|/|Z[\lambda]| \leq \alpha$, that is, a definite
proportion of parameters in $Z[\lambda]$ have escaping critical point).

\begin{rem} \label {compest}

Our estimates can be easily pushed further to obtain more precise
results.  For instance, it is clear that
$$
\alpha_{n+1} \leq \epsilon_{n+1}+(1-\epsilon_{n+1}) \epsilon_n
\sum_{k=0}^\infty \alpha_n^k \leq \epsilon_{n+1}+\frac {\epsilon_n}
{1-\alpha},
$$
so $\alpha_n \to 0$ (exponentially fast) for all parameters in $\NR^0$.
This in turn can be used to show that each
parameter in $\NR^0$ is a density point of the complement of $\MM$\footnote{
This last result does not hold for all parameters in $\NR$: there are
parameters $c \in \NR^\infty$ (well approximated by
cusps of little Mandelbrot sets) such that $\limsup |\D_\epsilon(c) \cap
\MM|/|\D_\epsilon(c)|=1$.
Our techniques show however that for every $c \in
\NR^\infty$, $\alpha_{n_k+1} \to 0$
(to see this, one needs to do the area estimate
jumping through central cascades, using a combinatorial procedure similar to
Theorem \ref {3.8}).
In particular, the upper density of the complement of the Mandelbrot set is
one at any $c \in \NR$: $\liminf |\D_\epsilon(c) \cap
\MM|/|\D_\epsilon(c)|=0$.
This result also follows from Graczyk-Swiatek's
proof of Shishikura's Theorem (personal communication by Jacek Graczyk).}.

\end{rem}


\end{document}